\documentclass[9pt]{article}

\usepackage{setspace}

\usepackage{amsmath,amssymb}
\usepackage{epsfig}
\usepackage{psfrag}

\title{Asymptotic expansion of the homogenized matrix in two weakly stochastic homogenization
settings}
\author{Ronan Costaouec\footnote{Ecole des Ponts, 77455 Marne-La-Vall\'ee Cedex 2 and
INRIA, MICMAC project, 78153 Le Chesnay Cedex, France, email: {\tt costaour@cermics.enpc.fr}}}
\date{\today}

\newcommand{\reels}{\mathbb{R}}

\newcommand{\esp}{\mathbb{E}}

\newcommand{\dom}{\mathcal{D}}

\def\og{\leavevmode\raise.3ex\hbox{$\scriptscriptstyle\langle\!\langle$~}}
\def\fg{\leavevmode\raise.3ex\hbox{~$\!\scriptscriptstyle\,\rangle\!\rangle$}}

\def\EE{\mathbb E}

\def\esp{\mathbb E}
\def\PP{\mathbb P}
\def\RR{\mathbb R}
\def\reels{\mathbb R}
\def\ZZ{\mathbb Z}

\def\longrightharpoonup{\relbar\joinrel\rightharpoonup}

\newtheorem{theo}{Theorem}[section]

\newtheorem{remarque}{Remark}[section]
\newtheorem{prop}{Proposition}[section]

\newtheorem{lemme}{Lemma}[section]
\newtheorem{corollaire}{Corollary}[section]
\newenvironment{preuve}{\textbf{Proof :}}{\begin{flushright}\rule{1ex}{1ex}\end{flushright}}
\newenvironment{preuve21}{\textbf{Proof of Proposition 2.1:}}{\begin{flushright}\rule{1ex}{1ex}\end{flushright}}
\newenvironment{preuvebis}{\textbf{Proof of Lemma 3.1:}}{\begin{flushright}\rule{1ex}{1ex}\end{flushright}}
\newenvironment{preuveter}{\textbf{Proof of Proposition 3.1:}}{\begin{flushright}\rule{1ex}{1ex}\end{flushright}}

\begin{document}

\maketitle

\begin{abstract}
This article studies some numerical approximations of the homogenized matrix for stochastic linear elliptic partial differential equations in divergence form: $-\mathrm{div} \left(A \left(\frac{\cdot}{\epsilon},\omega\right) \nabla u^\epsilon(\cdot , \omega) \right)=f$. We focus on the case when $A$ is a small perturbation $A_\eta$ of a reference periodic tensor $A_\mathrm{per}$, where $\eta$ encodes the size of the perturbation. In this case, it has been theoretically shown in~\cite{BlancLeBrisLions2,MPRF} for both models considered in this article that the exact homogenized matrix $A^\star_\eta$ possesses an expansion in powers of $\eta$, the coefficients of which are deterministic. In practice, one cannot manipulate the exact terms of such an expansion. All objects are subjected to a discretization approach for the variables $x$ (FEM) and $\omega$ (Monte Carlo method). Thus we need to derive a similar expansion for the \emph{approximated} random homogenized matrix. In contrast to the expansion of the exact homogenized matrix, the expansion of the approximate homogenized matrix contains intrinsically \emph{random} coefficients. In particular, the second order term is random in nature. The purpose of this work is to derive and study this expansion in function of the parameters of the approximation procedure (size of the truncated computational domain used, meshsize of the finite elements approximation).

\end{abstract}

\section{Introduction}
The context of this work is the homogenization of stochastic linear elliptic partial differential equations in divergence form 
\begin{eqnarray}
\label{eq:edp-epsilon}
\left\{
\begin{aligned}
-\mathrm{div} \left(A_\eta \left(\frac{\cdot}{\epsilon},\omega\right) \nabla u^\epsilon_\eta(\cdot , \omega) \right)=f \ \mathrm{on} \ \dom, \\
u^\epsilon_\eta(\cdot,\omega) = 0 \  \mathrm{on} \ \partial \dom,
\end{aligned}
\right.
\end{eqnarray}
where $\dom$ denotes a bounded domain of $\RR^d$, $d\in {\mathbb N}^*$ being the ambient dimension, and $f \in L^2(\mathcal{D})$. The matrix $A_\eta$ is random, symmetric, uniformly bounded and coercive, that is: 
\begin{equation}
\label{eq:hyp-bc-unif}
\left\{
\begin{array}{l}
\exists \gamma >0, \ \forall \xi \in \RR^d, \ \xi^T A_\eta(x, \omega) \xi \geq \gamma |\xi|^2 \ \mbox{almost everywhere and almost surely}, \\
\exists M >0 \mbox{ such that } \left\| A_\eta  \right\|_{\left(L^\infty(\RR^d \times \Omega)\right)^{d \times d}} \leq M,
\end{array}
\right.
\end{equation}
where $\gamma$ and $M$ do not depend on $\eta$.
Under additional stationarity hypotheses (the sense of which will be made precise below), it
is classical that, when $\epsilon \rightarrow 0$, the \emph{random} solutions $u^\epsilon_\eta$ of \eqref{eq:edp-epsilon} converge in some appropriate sense to a \emph{deterministic} function $u^\star_\eta$ solution to the following problem
\begin{eqnarray} \label{eq:edp-homog}
\left\{
\begin{aligned}
-\mathrm{div} \left(A^\star_\eta \nabla u^\star_\eta \right)=f \ \mathrm{on} \ \dom, \\
u^\star_\eta = 0 \  \mathrm{on} \ \partial \dom,
\end{aligned}
\right.
\end{eqnarray}
where $A^\star_\eta$ is the constant homogenized matrix. 
To compute $A^\star_\eta$, one needs first to compute some random functions, the so-called correctors $w_p^\eta$ (where $p \in \RR^d$), that are solutions to random auxiliary problems, namely the \emph{corrector problems} posed on $\RR^d$. As a consequence, solving these corrector problems, and hence computing $A^\star_\eta$, is in general computationally challenging. Appropriate truncations have to be considered, and the standard numerical strategies lead to extremely expensive computations. 
In this article, we do not deal with the very general framework of stochastic homogenization. We rather focus on the specific case when the matrix $A_\eta$ is a \emph{small} perturbation of an underlying periodic matrix denoted by $A_{\mathrm{per}}$. The parameter $\eta$ somehow encodes the degree of randomness of $A_\eta$: the smaller $\eta$, the less random $A_\eta$ is. When $\eta=0$, the matrix $A_\eta$ is deterministic. There are several ways to formalize the notion of a small perturbation in the context of stochastic homogenization (see~\cite{Anantharaman-Lebris-1,Anantharaman-Lebris-2,Singapour,BlancLeBrisLions1,BlancLeBrisLions2,MPRF,cras-ronan}). In this article, we will consider the following two models (for which all the details will be provided in Section 1.1 below):
\begin{itemize}
 \item \textbf{Model 1} [Standard discrete stationary setting for stochastic homogenization]: $A_\eta(x,\omega)$ is stationary, and the small perturbation is linear (see~\cite{MPRF}): $$A_\eta(x,\omega) = A_{\mathrm{per}}(x) + \eta A_1(x,\omega) + O(\eta^2),$$
where the expansion holds in some appropriate functional space.
 \item \textbf{Model 2} [Setting introduced in~\cite{BlancLeBrisLions1}]: the matrix reads $A_\eta = A_{\mathrm{per}} \left(\Phi_\eta^{-1}(x,\omega) \right)$, where $ \Phi_\eta$ is referred to as a stochastic diffeomorphism from $\RR^d$ to $\RR^d$ and $A_\mathrm{per}$ is a periodic matrix. An important point is that in this case $A_\eta$ is not necessarily stationary. The corresponding perturbative approach has been developed in~\cite{BlancLeBrisLions2}. The stochastic diffeomorphism $\Phi_\eta$ is supposed to be a linear perturbation of the identity, namely
$$\Phi_\eta(x,\omega) = x + \eta \Psi (x,\omega) + O(\eta^2),$$
in some appropriate functional space.
\end{itemize}

In both settings, it has been shown that the gradients of
the correctors $\nabla w_p^\eta$ and the homogenized matrix $A^\star_\eta$ possess an expansion in powers of $\eta$: $\nabla w_p^\eta=\nabla w_p^0 + \eta \nabla w_p^1 + O(\eta^2)$ and $A^\star_\eta = A^\star_{\mathrm{per}} + \eta A^\star_1 + O (\eta^2)$. Approximating $A^\star_\eta$ by $A^\star_{\mathrm{per}} + \eta A^\star_1$, one thus makes an error of the order of $\eta^2$. The functions $w_p^0$ and $w_p^1$ are solutions to  a deterministic and a random problem respectively. The dominant term $A^\star_{\mathrm{per}}$ in the expansion of $A^\star_\eta$ is nothing but the homogenized matrix associated to $A_{\mathrm{per}}$. The definition of $A^\star_1$ involves only $\nabla w_p^0$ and $  \nabla \EE \left(w^1_p \right)$. Since  $\EE \left(w^1_p \right)$ can be shown to solve a deterministic problem, the computation of the first order approximation $A^\star_{\mathrm{per}} + \eta A^\star_1$ only requires the solution of two \emph{deterministic} partial differential equations. This is far less demanding that solving a single \emph{stochastic} partial differential equation. This is the major advantage of these perturbative approaches. \\

Of course, in practice, one cannot exactly solve the corrector problems that are posed on an unbounded domain. As will be seen below, the classical discretization approach  consists in two steps: truncation of the corrector problem on a bounded computational domain of size $N$ and finite elements approximation with mesh size $h$. This introduces two types of error related to the parameters $N$ and $h$, respectively. Because of the truncation, the approximated homogenized matrix $A^{\star,h,N}_\eta$ is random. We will see that it also possesses an expansion in the variable $\eta$, almost surely. The main difference between this expansion and the exact one is that, at the discrete level, the difference $A^{\star,h,N}_\eta - A^{\star,h,N}_\mathrm{per}-\eta A^{\star,h,N}_1$ is random in nature, and so is the second order error estimate. Therefore, in order to ensure the relevance of the first order approximation of  $A^{\star,h,N}_\eta$ , one needs to bound the random error $\eta^{-2} (A^{\star,h,N}_\eta - A^{\star,h,N}_\mathrm{per}-\eta A^{\star,h,N}_1 )$ of second order, in some appropriate probability space, namely $\left(L^\infty(\Omega)\right)^{d \times d}$ in the sequel. In the setting of Model 2, this question has already been addressed from a numerical point of view in~\cite{cras-ronan}.  The aim of the present work is to theoretically derive rigorous bounds on the second order error for both models and to understand how these bounds depend on the parameters of the discretization procedure. To this end, we first have to make precise the appropriate functional spaces in which the original expansions
%% of the original matrix $A_\eta$  and of the stochastic diffeomorphism $\Phi_\eta$ 
hold. We next deduce for each model in which sense the second order error terms are bounded. Our presentation elaborates on the previous works~\cite{BlancLeBrisLions1,BlancLeBrisLions2}. We only provide the formalism for self-consistency and we refer to~\cite{BlancLeBrisLions2} for more details.

\subsection{Probabilistic setting}

Throughout this article, $(\Omega, {\mathcal F}, \PP)$ is a
probability space and we denote by $\EE(X) =
\int_\Omega X(\omega) d\PP(\omega)$ the expectation value of any
random variable $X\in L^1(\Omega, d\PP)$. We assume that the group $(\ZZ^d, +)$ acts on
$\Omega$. We denote by $(\tau_k)_{k\in \ZZ^d}$ this action, and
assume that it preserves the measure $\PP$, that is, for all
$\displaystyle  k\in \ZZ^d$ and all $A \in {\cal F}$, $\displaystyle \PP(\tau_k A)
  = \PP(A)$. We assume that the action $\tau$ is {\em ergodic}, that is,
  if $A \in {\mathcal F}$ is such that $\tau_k A = A$ for any $k \in
  \ZZ^d$, then $\PP(A) = 0$ or 1. 
In addition, we define the following discrete notion of stationarity (see \cite{BlancLeBrisLions1}): any
$F\in L^1_{\rm loc}\left(\RR^d, L^1(\Omega)\right)$ is said to be
{\em stationary} if, for all $k\in \ZZ^d$,
\begin{equation}
\label{eq:stationnarite-disc}
F(x+k, \omega) = F(x,\tau_k\omega), \ \mbox{almost everywhere and almost surely.}
\end{equation}

Note that this setting is a straightforward classical variant of the more commonly used (continuous) stationary setting for random homogenization, for which the shift $\tau$ is indexed by elements of the group $\RR^d$ and \eqref{eq:stationnarite-disc} holds for all $k \in \RR^d$ instead of $\ZZ^d$.
In our discrete setting, the ergodic theorem~\cite{krengel,shiryaev} can be
stated as follows: 

\begin{theo}
\label{theo:ergodic}
Let $F\in L^\infty\left(\RR^d, L^1(\Omega)\right)$ be a stationary random
variable in the above sense. For $k = (k_1,k_2, \dots k_d) \in \ZZ^d$,
we set $\displaystyle |k|_\infty = \sup_{1\leq i \leq d} |k_i|$. Then
$$
\frac{1}{(2N+1)^d} \sum_{|k|_\infty\leq N} F(x,\tau_k\omega)
\mathop{\longrightarrow}_{N\rightarrow \infty}
\EE\left(F(x,\cdot)\right) \quad \mbox{in }L^\infty(\RR^d),
\mbox{ almost surely}.
$$
This implies that (denoting by $Q$  the unit cube in $\RR^d$)
$$
F\left(\frac x \varepsilon ,\omega \right)
\mathop{\longrightharpoonup}_{\varepsilon \rightarrow 0}^*
\EE\left(\int_Q F(x,\cdot)dx\right) \quad \mbox{in }L^\infty(\RR^d),
\mbox{ almost surely}.
$$
\end{theo}

\subsection{Stochastic homogenization results}
Standard results of stochastic homogenization (see~\cite{blp,jikov}) apply to problem~\eqref{eq:edp-epsilon} when the matrix $A_\eta$ is stationary. For any fixed value of $\eta$, they provide the following result. 
\begin{theo}
We consider Model 1. Suppose that $A_\eta$ is a symmetric matrix, stationary in the sense of~\eqref{eq:stationnarite-disc}, uniformly bounded and coercive in the sense of~\eqref{eq:hyp-bc-unif}. Then the homogenized matrix $A^\star_\eta$ appearing in \eqref{eq:edp-homog}
is defined by
\begin{equation}\label{eq:Cas1-Astar}
\left[ A^\star_\eta \right]_{ij} = \EE \left[ \int_{Q} e_i^T A_\eta \left(e_j + \nabla w^\eta_{e_j} \right)  \right],
\end{equation}
where, for any $p \in \RR^d$, $w^\eta_p$ denotes the unique solution (up to the addition of a random constant) of the corrector problem
\begin{equation} \label{eq:Cas1-PC-Linear}
\left\{
\begin{array}{l}
\displaystyle{-\mathrm{div} \left(A_\eta  \left(p+ \nabla w^\eta_p \right) \right)=0 \ \mathrm{on} \ \RR^d, \ \mbox{almost surely},} \\
\displaystyle{\nabla w_p^\eta\ \ \mbox{is stationary in the sense of~\eqref{eq:stationnarite-disc} },} \\
\displaystyle{\EE \left[ \int_Q \nabla w_p^\eta \right] = 0.}
\end{array}
\right.
\end{equation}
\end{theo}
In the sequel we will always assume that $\left| p\right|=1$.\\

In the case $A_\eta(x,\omega) = A_{\mathrm{per}} \left(\Phi_\eta^{-1}(x,\omega) \right)$, there exists an analogous result due to Blanc, Le Bris and Lions (see~\cite{BlancLeBrisLions1}). Its statement requires to make precise the notion of \emph{stochastic diffeomorphism} mentioned above.  The map $\Phi_\eta$ is said to be a stochastic
diffeomorphism if it satisfies
\begin{eqnarray}
\label{eq:hyp-diffeo-as}
&&\Phi_\eta(\cdot,\omega) \ \mbox{is a diffeomorphism almost surely},\\
\label{eq:hyp-diffeo-1}
&&\underset{x \in \reels^d, \omega \in \Omega}{\mathrm{Ess}\inf} \left( \mathrm{det} \left(\nabla \Phi_\eta (x,\omega) \right)\right)=\nu > 0, \\
\label{eq:hyp-diffeo-2}
&&\underset{x \in \reels^d, \omega \in \Omega}{\mathrm{Ess}\sup} \left( |\nabla \Phi_\eta(x,\omega)| \right) = M' < \infty, \label{HDiff2}\\
\label{eq:hyp-diffeo-3}
&&\nabla \Phi_\eta \ \mbox{is stationary in the sense of~\eqref{eq:stationnarite-disc}}. \label{HDiff3}
\end{eqnarray}
Under these hypotheses on $\Phi_\eta$, the following theorem gives the homogenized problem associated to~\eqref{eq:edp-epsilon}.
\begin{theo}
We consider Model 2. Suppose that $A_\eta=A_{\mathrm{per}} \circ \Phi_\eta^{-1}$ where $A_\mathrm{per}$ is a symmetric matrix, uniformly bounded and coercive, and $\Phi_\eta$ satisfies~\eqref{eq:hyp-diffeo-as},~\eqref{eq:hyp-diffeo-1},~\eqref{eq:hyp-diffeo-2} and~\eqref{eq:hyp-diffeo-3}.  Then the homogenized matrix $A^\star_\eta $ appearing in \eqref{eq:edp-homog} is defined by
\begin{equation}\label{eq:Cas2-Astar}
\left[ A^\star_\eta \right]_{ij} =\left(\EE \left[ \int_{Q} \mathrm{det} \left(\nabla \Phi_\eta \right) \right] \right)^{-1} \EE \left[ \int_{Q} \mathrm{det} \left( \nabla \Phi_\eta \right) e_i^T A_\mathrm{per} \left(e_j + \left(\nabla \Phi_\eta \right)^{-1} \nabla w^\eta_{e_j} \right)  \right],
\end{equation}
where, for any $p \in \RR^d$, $w^\eta_p$ denotes the unique solution (up to the addition of a random constant) of the corrector problem
\begin{equation} \label{eq:Cas2-PC-Linear}
\left\{
\begin{array}{l}
\displaystyle{-\mathrm{div} \left( \mathrm{det} \left( \nabla \Phi_\eta (\cdot,\omega) \right)\left(\nabla \Phi_\eta (\cdot,\omega) \right)^{-T} A_\mathrm{per} \left(p+ \left(\nabla \Phi_\eta (\cdot,\omega) \right)^{-1}\nabla w^\eta_p (\cdot,\omega) \right) \right)=0 \ \mathrm{on} \ \RR^d, \ \mbox{almost surely},} \\
\displaystyle{\nabla w_p^\eta\ \ \mbox{is stationary in the sense of~\eqref{eq:stationnarite-disc}},} \\
\displaystyle{\EE \left[ \int_Q \nabla w_p^\eta \right] = 0.}
\end{array}
\right.
\end{equation}
\end{theo}

\subsection{Standard numerical approximation}

In practice, problems~\eqref{eq:Cas1-PC-Linear} and~\eqref{eq:Cas2-PC-Linear} are solved numerically.
The first step is to introduce a truncation. Following~\cite{BourgeatPiatnitski}, we approximate~\eqref{eq:Cas1-PC-Linear} by
\begin{equation} \label{eq:Cas1-PC-Linear-N}
\left\{
\begin{array}{l}
\displaystyle{-\mathrm{div} \left(A_\eta (\cdot,\omega)  \left(p+ \nabla w^{\eta,N}_p(\cdot,\omega) \right) \right)=0 \ \mathrm{on} \ Q_N \ \mbox{almost surely},} \\
\displaystyle{ w_{p}^{\eta,N}(\cdot,\omega)\ \ \mbox{is} \ Q_N-\mbox{periodic},}
\end{array}
\right.
\end{equation} 
where $Q_N = \left[-N-\frac{1}{2}, N + \frac{1}{2} \right]^d$. We consider an analogous truncated problem for~\eqref{eq:Cas2-PC-Linear}. A classical finite element method is then used to approximate the solutions of~\eqref{eq:Cas1-PC-Linear-N}. We consider a periodic triangulation
$\mathcal{T}^{(Q)}_h$ of the unit cell $Q=\left[-\frac{1}{2},\frac{1}{2}\right]^d$. Replicating it, we obtain
a triangulation 
$$\mathcal{T}_h = \cup_{k \in \ZZ^d} \left( k + \mathcal{T}^{(Q)}_h
\right)$$ of $\RR^d$. We denote by $V_h^{\rm per}(Q_N)$ the space of functions $\varphi_h$ defined on
$\RR^d$, $Q_N$-periodic, whose restriction to $Q_N$ is in a $\mathbb{P}_1$-Lagrange finite element space built from
$\mathcal{T}_h^N = \mathcal{T}_h \cap Q_N$, and which satisfy
\begin{equation}
\int_{Q_N} \varphi_h = 0. \nonumber
\end{equation}
Let $\left\lbrace \phi_k \right\rbrace_{1 \leq k \leq N_v} $ be a basis of $V^h_{\mathrm{per}} \left(Q_N \right)$:
$$
V^h_{\mathrm{per}} \left(Q_N \right)= \mbox{span} \left(\left\lbrace \phi_k \right\rbrace_{1 \leq k \leq N_v}\right),
$$
where $N_v=N_v(N)$ is the number of degrees of freedom considered.\\

In the standard case (Model 1), we define the \emph{approximated corrector} $w^{\eta,h,N}_p$ as the solution to the variational formulation
\begin{equation} 
\label{eq:correcteur-hN}
\left\{
\begin{array}{l}
\mbox{Find} \ 
w_p^{\eta,h,N}(\cdot,\omega) \in V_h^{\rm per}(Q_N) \ \ 
\mbox{such that,} \ 
\\
\displaystyle{\forall \varphi_h \in V^{\rm per}_h(Q_N), \ \ \int_{Q_N}
  A_\eta(\cdot,\omega)
\left( p +  \nabla w_p^{\eta,h,N}(\cdot,\omega)\right) \cdot \nabla \varphi_h = 0 \ \ \mbox{almost surely,}} 
\end{array}
\right.
\end{equation}
and the \emph{approximated homogenized matrix} by
\begin{equation}
\label{eq:homog-matrix-Nh}
\forall 1 \leq i,j \leq d, \ \ \left[A^{\star,h,N}_\eta \right]_{ij} (\omega) =\frac{1}{|Q_N|} \int_{Q_N} e_i^T A_\eta (\cdot,\omega) \left(e_j + \nabla w_{e_j}^{\eta,h,N}(\cdot,\omega) \right).
\end{equation}
In the second setting (Model 2), we similarly define the approximated corrector as the solution to
\begin{equation} 
\label{eq:correcteur-hN-cas2}
\left\{
\begin{array}{l}
\mbox{Find} \ 
w_p^{\eta,h,N}(\cdot,\omega) \in V_h^{\rm per}(Q_N) \ \ 
\mbox{such that, for any } \varphi_h \in V^{\rm per}_h(Q_N),\ 
\\
\displaystyle{\int_{Q_N}
  \mathrm{det}(\nabla \Phi_\eta(\cdot,\omega))(\nabla \Phi_\eta(\cdot,\omega))^{-T} A_{\mathrm{per}}
\left( p + (\nabla \Phi_\eta(\cdot,\omega))^{-1} \nabla w_p^{\eta,h,N}(\cdot,\omega)\right) \cdot \nabla \varphi_h = 0 \ \ \mbox{almost surely}},
\end{array}
\right.
\end{equation}
and the homogenized matrix is defined by
\begin{equation}
\label{eq:a-hN-cas2}
\begin{array}{l}
\- \left[A^{*,h,N}_\eta \right]_{ij}(\omega) 
= \mathrm{det} \left( \int_{Q_N} \nabla \Phi_\eta(\cdot,\omega)  \right)^{-1} \left(  \int_{Q_N} \mathrm{det}\left(\nabla \Phi_\eta(\cdot,\omega) \right)
\left(e_i + \left(\nabla \Phi_\eta(\cdot,\omega) \right)^{-1} \nabla w^{\eta,h,N}_{e_i}(\cdot,\omega) \right)^T A_{\mathrm{per}} e_j  \right).
\end{array}
\end{equation}
Note that for both~\eqref{eq:correcteur-hN} and~\eqref{eq:correcteur-hN-cas2} the solution $w_p^{\eta,h,N}$ is a \emph{random} field because of the truncation procedure. It follows that both~\eqref{eq:homog-matrix-Nh} and~\eqref{eq:a-hN-cas2} are \emph{random} matrices. It is only in the limit $N=\infty$ that these objects become \emph{deterministic}. In the numerical practice, one then commonly considers that the best approximation of $A^\star_\eta$ is given by $\EE \left(A^{*,h,N}_\eta \right)$ which is in turn estimated using an empirical mean, computed with standard Monte-Carlo methods.

\section{Expansion of the random matrix $A_\eta$ (Model 1)}

In this section, we consider the weakly stochastic setting of Model 1. We will turn our attention to Model 2 in Section 3. Our main result for Model 1 is  Proposition~\ref{prop:cas1-dl-hN}. As announced above, it makes precise the behaviour of the \emph{random} second order error in the expansion of the approximated homogenized matrix $A^{\star,h,N}_\eta$, under appropriate hypotheses. Then, passing to the limit $h \rightarrow 0$, we prove that this result extends to the approximated homogenized matrix when only truncation is taken into account. Finally, letting $N$ go to infinity, we recover the expansion of the exact homogenized matrix $A^\star_\eta$ derived in~\cite{MPRF}. The functional setting described below is simple and the arguments we use in the proof of Proposition 2.1 are standard.
Our aim is to illustrate in a simple and relevant framework that the randomness of the second order error in the expansions of $\nabla w^{\eta,h,N}_p$ and $A^{\star,h,N}_\eta$ does not affect the validity of the approximation $A^{\star,h}_\mathrm{per} + \eta \EE (A_1^{\star,h,N})$. In particular, this amounts to prove that the quantity $\eta^{-2} (A^{\star,h,N}_\eta - A^{\star,h}_\mathrm{per}-\eta A^{\star,h,N}_1 )$ is bounded independently of $h$, $N$ and $\eta$ in some probability space. This is precisely what ensures Proposition 2.1 below. 

\subsection{Assumptions}

We suppose that the matrix $A_\eta$ admits in $\left(L^\infty ( Q \times \Omega)\right)^{d \times d}$ the expansion
\begin{equation}
\label{eq:dl}
A_\eta(x,\omega) = A_\mathrm{per}(x) + \eta A_1(x,\omega) + R_\eta (x,\omega),
\end{equation}
where $A_\mathrm{per} \in \left(L^\infty(Q) \right)^{d \times d}$, $A_1 \in \left( L^\infty (Q \times \Omega) \right)^{d \times d}$ and $R_\eta = O(\eta^2)$ in $\left( L^\infty (Q \times \Omega) \right)^{d \times d} $.
This means that  
\begin{eqnarray}
\label{eq:dl-H1}
\underset{\eta \rightarrow 0}{\lim} \left\|A_\eta - A_\mathrm{per} \right\|_{\left( L^\infty (Q \times \Omega) \right)^{d \times d}}&=&0,  \\
\label{eq:dl-H2}
\underset{\eta \rightarrow 0}{\lim} \left\| \eta^{-1}\left(A_\eta - A_\mathrm{per}\right) -A_1 \right\|_{\left( L^\infty (Q \times \Omega) \right)^{d \times d}}&=&0,
\end{eqnarray}
and there exists a \emph{deterministic} constant $ C_R $ independent of $\eta$ such that, when $| \eta | \leq 1 $, 
\begin{equation}
\label{eq:dl-H3}
\eta^{-2}\left\|R_\eta \right\|_{\left( L^\infty (Q \times \Omega)\right)^{d \times d} } \leq C_R, 
\end{equation}
where we recall that 
\begin{equation}
\left\| v \right\|_{ \left(L^\infty (Q \times \Omega)\right)^{d \times d}} = \underset{(x,\omega) \in Q \times \Omega}{\mathrm{Ess}\sup} \left|v(x,\omega )\right|. \nonumber
\end{equation}
From~\eqref{eq:hyp-bc-unif}, we know that $A_\eta$ is bounded and uniformly coercive. 
Using \eqref{eq:dl-H1} and~\eqref{eq:hyp-bc-unif}, we see that 
\begin{equation}
\label{eq:hyp-b-unif-0}
\left\|A_\mathrm{per}  \right\|_{\left(L^\infty (Q) \right)^{d \times d}}=\underset{\eta \rightarrow 0}{\lim} \left\| A_\eta \right\|_{\left(L^\infty (Q \times \Omega) \right)^{d \times d}} \leq M , 
\end{equation}
and, for any $\xi \in \RR^d$, we have
\begin{equation}
\label{eq:hyp-c-unif-0}
\xi^T A_\mathrm{per}(x) \xi =  \underset{\eta \rightarrow 0}{\lim} \xi^T A_\eta (x,\omega) \xi \geq \gamma \left|\xi \right|^2 \ \mbox{almost everywhere}. 
\end{equation}
The matrix $A_\mathrm{per}$ is thus also uniformly bounded and coercive. \\

Assuming that the expansion of $A_\eta$ holds in $\left(L^\infty \left(Q \times \Omega \right) \right)^{d \times d}$ is 
relevant from the point of view of modelization, and it somehow simplifies the proof of Proposition 2.1. As will be seen below (Remark 2.1), up to slight modifications, our main result (Proposition 2.1) extends to the case when the original expansion of $A_\eta$ holds in a weaker sense, namely when it holds almost surely with uniformly integrable bounds. Of course, in this case, the second order error in the expansion of $A^{\star,h,N}_\eta$ is bounded in a weaker sense. Lastly, it is to be mentioned that our proof does not directly extend to the case when the expansion of $A_\eta$ holds in spaces of the form $\left(L^\infty \left(Q; L^p( \Omega )\right) \right)^{d \times d}$ with $ 0 < p \leq + \infty$; which models the idea of possibly large but rare local perturbations (see~\cite{Singapour} for a detailed presentation of this latter model).

\subsection{Formal expansion}

Following the method introduced in~\cite{BlancLeBrisLions2}, we first postulate the formal expansion of the solution to~\eqref{eq:correcteur-hN}:
\begin{equation}
\label{eq:dl-w}
w_p^{\eta,h,N} = w^{0,h,N}_p + \eta w^{1,h,N}_p + r_p^{\eta,h,N},
\end{equation}
where $\nabla r_p^{\eta,h,N} = O (\eta^2)$ in some appropriate space. We will sucessively identify $w^{0,h,N}_p$ and $w^{1,h,N}_p$ and prove the validity of the expansion. Formally inserting this expression in~\eqref{eq:correcteur-hN}, we obtain that the function $w_p^{0,h,N}$ is independent of $N$ (we denote it by $w^{0,h}_p$ in the sequel) and solves
\begin{equation} 
\label{eq:correcteur-0-hN}
\left\{
\begin{array}{l}
\mbox{Find} \ 
w_p^{0,h} \in V_h^{\rm per}(Q) \ \ 
\mbox{such that,} \ 
\\
\displaystyle{\forall \varphi_h \in V^{\rm per}_h(Q), \ \ \int_{Q}
  A_{\mathrm{per}}
\left( p +  \nabla w_p^{0,h}\right) \cdot \nabla \varphi_h = 0,
}
\end{array}
\right.
\end{equation}
and that the function $w_{p}^{1,h,N}$ is solution to
\begin{equation} 
\label{eq:correcteur-1-hN}
\left\{
\begin{array}{l}
\mbox{Find} \ 
w_p^{1,h,N}(\cdot,\omega) \in V_h^{\rm per}(Q_N) \ \ 
\mbox{such that,} \ 
\\
\displaystyle{\forall \varphi_h \in V^{\rm per}_h(Q_N), \ \ \int_{Q_N} A_{\mathrm{per}} \nabla w_p^{1,h,N}(\cdot,\omega) \cdot \nabla \varphi_h + \int_{Q_N}
  A_1(\cdot,\omega)
\left( p +  \nabla w_p^{0,h}\right) \cdot \nabla \varphi_h = 0 \ \ \mbox{almost surely.}
}
\end{array}
\right.
\end{equation}
In addition, substituting the expansions of $A_\eta$ and $w_p^{\eta,h,N}$ into~\eqref{eq:homog-matrix-Nh}, we formally obtain  
\begin{eqnarray}
\label{eq:Dev-Astar-1}
A^{\star,h,N}_\eta (\omega) = A^{\star,h}_\mathrm{per} + \eta A_1^{\star,h,N}(\omega) + O(\eta^2),
\end{eqnarray}
where the terms of order zero and one are respectively defined by
\begin{eqnarray}
\label{eq:a-0-hN}
& &\left[ A_\mathrm{per}^{\star,h}\right] _{ij} = \int_Q e_i^T A_\mathrm{per} \,
\left( e_j + \nabla w_{e_j}^{0,h} \right), \\
\label{eq:a-1-hN}
& &\left[ A^{\star,h,N}_1\right] _{ij}  (\omega) = \frac{1}{|Q_N|} \int_{Q_N} e_i^T A_\mathrm{per} \nabla w^{1,h,N}_{e_j} (\cdot,\omega)
+ \frac{1}{|Q_N|} \int_{Q_N} e_i^T A_1 (\cdot,\omega) \left( \nabla w^{0,h}_{e_j} + e_j \right).
\end{eqnarray}
In the sequel, we will make precise and rigorously justify the expansions~\eqref{eq:dl-w} and~\eqref{eq:Dev-Astar-1}.

\subsection{Main result}

Our main result in this section is the following.
\begin{prop}
\label{prop:cas1-dl-hN}
Suppose that $A_\eta$ is a symmetric matrix that satisfies (\ref{eq:hyp-bc-unif}) and is stationary in the sense of~\eqref{eq:stationnarite-disc}. Suppose, in addition, that it satisfies~(\ref{eq:dl}), (\ref{eq:dl-H1}), (\ref{eq:dl-H2}). We assume that (\ref{eq:dl-H3}) holds, namely the second order error is $O(\eta^2)$ in $\left(L^\infty (Q \times \Omega) \right)^{d \times d}$. Then there exists a constant $C$ independent of $\eta$, $\omega$, $N$ and $h$, such that, for $|\eta| \leq 1$,
\begin{equation}
\label{eq:cas1-dl-w-hN}
 \eta^{-2} \left\| 
\nabla w^{\eta,h,N}_p(\cdot,\omega) - \nabla w^{0,h}_p - \eta \nabla
w^{1,h,N}_p(\cdot,\omega) \right\|_{\left(L^2(Q_N)\right)^{d}} \leq C \sqrt{| Q_N |} \ \ \mbox{almost surely},
\end{equation}
where $w^{\eta,h,N}_p$, $w^{0,h}_p$ and $w^{1,h,N}_p$ are solutions
to~(\ref{eq:correcteur-hN}), (\ref{eq:correcteur-0-hN})
and~(\ref{eq:correcteur-1-hN}), respectively, and such that
\begin{equation}
\label{eq:cas1-dl-a-hN}
 \eta^{-2} 
\left| A^{\star,h,N}_{\eta}(\omega) - A_\mathrm{per}^{\star,h} - \eta A^{\star,h,N}_1(\omega) \right|
\leq C \ \ \mbox{almost surely},
\end{equation}
where $A^{\star,h,N}_{\eta}$, $A_\mathrm{per}^{\star,h}$ and $A^{\star,h,N}_1$ are defined by~\eqref{eq:homog-matrix-Nh},~\eqref{eq:a-0-hN} and~\eqref{eq:a-1-hN}.
\end{prop}

Note that, as the constant $C$ in~\eqref{eq:cas1-dl-w-hN} and~\eqref{eq:cas1-dl-a-hN} is independent
of $\omega$, the expansions of  $\nabla w_p^{\eta,h,N}$ and $A^{\star,h,N}_\eta$ hold in $\left(L^\infty \left(\Omega;L^2(Q_N) \right)\right)^{d}$ and $\left( L^\infty (\Omega)\right)^{d \times d}$, respectively.\\
%\end{remarque}

\begin{preuve21}
Our goal is to prove estimates \eqref{eq:cas1-dl-w-hN} and~\eqref{eq:cas1-dl-a-hN}. To do so, following the methodology of~\cite{BlancLeBrisLions2}, we begin with justifying the expansion~\eqref{eq:dl-w} in our discrete framework. This is the purpose of Steps 1, 2 and 3. First we will check that $\nabla w^{\eta,h,N}_p$ is bounded independently of $\eta,\omega,N$ and $h$ and that it converges to the gradient $\nabla w^{0,h}_p$ of a \emph{deterministic} function (Step 1). We will then verify that the
first order error term of $\nabla w^{\eta,h,N}_p$, namely $\nabla v^{\eta,h,N}_p$ defined below, is bounded independently of $\eta,\omega,N$ and $h$ and that it converges to the gradient $\nabla w^{1,h,N}_p$ of a \emph{random} function (Step 2). In Step 3, we prove that the second order error of $\nabla w^{\eta,h,N}_p$, namely $\nabla z^{\eta,h,N}_p$ defined below, is bounded independently of $\eta,\omega,N$ and $h$. In other terms we prove \eqref{eq:cas1-dl-w-hN}, using bounds derived at Steps 1 and 2. In Step 4, remarking that the second order error term of the homogenized matrix $A^{\star,h,N}_\eta$ depends on $\nabla w^{\eta,h,N}_p$, $\nabla v^{\eta,h,N}_p$ and $\nabla z^{\eta,h,N}_p$, we use the bounds from Steps 1, 2 and 3 to prove~\eqref{eq:cas1-dl-a-hN}.

\paragraph{\underbar{Step 1:}} Our goal is first to prove that $\nabla w^{\eta,h,N}_p$ is bounded in $\left(L^2(Q_N)\right)^{d}$, almost surely and independently of $\eta,\omega,N$ and $h$, and next to show that $w^{\eta,h,N}_p$ converges to a \emph{deterministic} function $w^{0,h}_p$ independent of $N$. Choosing $\varphi_h = w^{\eta,h,N}_p$ as test function in~\eqref{eq:correcteur-hN} and using~\eqref{eq:hyp-bc-unif},
we obtain that 
\begin{equation}
\label{eq:W}
\gamma \left\| \nabla w^{\eta,h,N}_p(\cdot,\omega) \right\|^2_{\left(L^2 (Q_N)\right)^{d}} \leq \left(\int_{Q_N} | A_\eta(\cdot,\omega) p |^2 \right)^{\frac{1}{2}} \left(\int_{Q_N} | \nabla w^{\eta,h,N}_p |^2 \right)^{\frac{1}{2}}, \nonumber
\end{equation}
which implies that (recall that $|p|=1$)
\begin{equation}
\label{eq:bound-w}
\left\| \nabla w^{\eta,h,N}_p (\cdot,\omega) \right\|_{\left(L^2 (Q_N)\right)^{d}}  \leq \frac{\left\|A_\eta \right\|_{\left(L^\infty(Q \times \Omega)\right)^{d \times d}}}{\gamma} |Q_N|^{\frac{1}{2}}. \nonumber
\end{equation}
Using~(\ref{eq:hyp-bc-unif}), we deduce that
\begin{equation}
\label{eq:Cw}
\left\| \nabla w^{\eta,h,N}_p (\cdot,\omega) \right\|_{\left(L^2 (Q_N)\right)^d} \leq  C_w |Q_N|^{\frac{1}{2}},
\end{equation}
where 
\begin{equation}
\label{eq:def-Cw}
C_w= \frac{M}{\gamma} \nonumber
\end{equation}
is a constant independent of $\eta,\omega,N$ and $h$.
Thus $ \nabla w^{\eta,h,N}_p (\cdot,\omega)$ is bounded uniformly in $\left(L^2 (Q_N)\right)^{d}$ independently from $\eta$ and almost surely.
We then remark that
$$ \nabla w^{\eta,h,N}_p (\cdot,\omega) = \sum_{i=1}^{N_v} \alpha^\eta_{p,i}(\omega) \nabla \phi_i,$$
where $\alpha^\eta_p = (\alpha^\eta_{p,i})_{1 \leq i \leq N_v}$ denotes the coordinates of $ w^{\eta,h,N}_p (\cdot,\omega)$ in the basis of $V_h^{\mathrm{per}}(Q_N)$. In the finite dimensional space $V^\mathrm{per}_h(Q_N)$, all norms are equivalent and we deduce from~\eqref{eq:Cw} that $\alpha^\eta_p$ is bounded independently of $\eta$ in $\RR^{N_v}$. Thus, up to extracting a subsequence, $\alpha^\eta_p$ converges to a vector $\alpha^0_p$ in $\RR^{N_v}$ almost surely. Consequently, we see that
$$ \nabla w^{\eta,h,N}_p (\cdot,\omega) \mathop{\longrightarrow}_{\eta \rightarrow 0} \ \sum_{i=1}^{N_v} \alpha^0_{p,i}(\omega) \nabla \phi_i = \nabla \left( \sum_{i=1}^{N_v} \alpha^0_{p,i}(\omega)\phi_i\right)  = \nabla w_p^{0,h,N} (\cdot,\omega) \ \mbox{in} \left(L^2(Q_N)\right)^{d} \ \mbox{almost surely,}
$$
where $ w^{0,h,N}_p(\cdot,\omega) \in V^{\mathrm{per}}_h(Q_N)$ is defined by 
$$w^{0,h,N}_p(\cdot,\omega):= \sum_{i=1}^{N_v} \alpha^0_{p,i}(\omega)\phi_i .$$
We now return to~(\ref{eq:correcteur-hN}) which we decompose using the expansion~\eqref{eq:dl} of $A_\eta$: for any $\varphi_h \in V^{\mathrm{per}}_h (Q_N), $ 
\begin{equation}
\label{eq:diff-0-hN}
 \ \int_{Q_N} A_\mathrm{per} \left(p+ \nabla w_p^{\eta,h,N}(\cdot , \omega) \right) \cdot \nabla \varphi_h + \eta \int_{Q_N} A_1(\cdot,\omega) \left(p+ \nabla w_p^{\eta,h,N}(\cdot , \omega) \right) \nabla \varphi_h + \int_{Q_N} R_\eta(\cdot,\omega)  \left(p+ \nabla w_p^{\eta,h,N}(\cdot , \omega) \right)\cdot \nabla \varphi_h =0.
\end{equation}
We are going to pass to the limit $\eta \rightarrow0$ in~\eqref{eq:diff-0-hN}. By definition of $\nabla w^{0,h,N}_p(\cdot,\omega)$, 
\begin{equation}
\int_{Q_N} A_\mathrm{per} \left(p+ \nabla w_p^{\eta,h,N}(\cdot , \omega) \right) \cdot \nabla \varphi_h \underset{\eta \rightarrow 0}{\longrightarrow} \int_{Q_N} A_\mathrm{per} \left(p+ \nabla w_p^{0,h,N}(\cdot , \omega) \right) \cdot \nabla \varphi_h \nonumber
\end{equation}
almost surely, and
\begin{equation}
\int_{Q_N} A_1(\cdot,\omega) \left(p+ \nabla w_p^{\eta,h,N}(\cdot , \omega) \right) \cdot \nabla \varphi_h \underset{\eta \rightarrow 0}{\longrightarrow} \int_{Q_N} A_1(\cdot,\omega) \left(p+ \nabla w_p^{0,h,N}(\cdot , \omega) \right) \cdot \nabla \varphi_h \nonumber
\end{equation}
almost surely. Using \eqref{eq:dl-H3}, we see that, when $|\eta| \leq 1$,
\begin{eqnarray}
\left| \int_{Q_N} R_\eta(\cdot,\omega)  \left(p+ \nabla w_p^{\eta,h,N}(\cdot , \omega) \right)\cdot \nabla \varphi_h \right| &\leq&  \left\|R_\eta \right\|_{\left(L^\infty (Q \times \Omega)\right)^{d \times d}} \left\|p + \nabla w_p^{\eta,h,N}(\cdot,\omega)\right\|_{\left(L^2(Q_N)\right)^d} \left\| \nabla \varphi_h \right\|_{\left(L^2(Q_N)\right)^d} \nonumber \\
&\leq& C_R \eta^2 \left\|p + \nabla w_p^{\eta,h,N}(\cdot,\omega)\right\|_{\left(L^2(Q_N)\right)^d} \left\| \nabla \varphi_h \right\|_{\left(L^2(Q_N)\right)^d} \nonumber \\
&\leq& C_R \eta^2 \left( \left\|p\right\|_{\left(L^2(Q_N)\right)^d} + C_w |Q_N|^{\frac{1}{2}} \right) \left\| \nabla \varphi_h \right\|_{\left(L^2(Q_N)\right)^d}, \nonumber
\end{eqnarray}
where we recall that $C_w$ is independent of $\eta$.
We thus obtain that
$$\underset{\eta \rightarrow 0}{\lim}  \int_{Q_N} R_\eta  \left(p+ \nabla w_p^{\eta,h,N}(\cdot , \omega) \right)\cdot \nabla \varphi_h  = 0.$$
Passing to the limit $\eta \rightarrow 0$ in~\eqref{eq:diff-0-hN}, we obtain that $w^{0,h,N}_p(\cdot,\omega)$ satisfies 
\begin{equation}
\label{eq:correcteur-0-hN-bis}
\displaystyle{ \forall \varphi_h \in V_h^{\mathrm{per}} (Q_N), \ \ \int_{Q_N} A_{\mathrm{per}} \nabla w^{0,h,N}_p(\cdot,\omega) \cdot \nabla \varphi_h+ \int_{Q_N} A_{\mathrm{per}} p \cdot \nabla \varphi_h = 0. } 
\end{equation}
Note that this equation has a unique solution in $V^h_\mathrm{per}(Q_N)$. We now show that $w^{0,h,N}_p = w^{0,h}_p$, where $w^{0,h}_p$ is the solution of~\eqref{eq:correcteur-0-hN}. Indeed, for any $\varphi_h \in V_h^{\mathrm{per}} (Q_N) $, we see that
\begin{eqnarray}
\int_{Q_N}  A_{\mathrm{per}} \left( \nabla w^{0,h}_p+ p \right) \cdot \nabla \varphi_h &=& \sum_{|k|_{\infty} \leq N} \int_{Q +k}  A_{\mathrm{per}} \left( \nabla w^{0,h}_p + p \right)\cdot \nabla \varphi_h \nonumber \\
&=& \sum_{|k|_{\infty} \leq N} \int_{Q}  A_{\mathrm{per}}  \left(\nabla w^{0,h}_p+ p \right) \cdot \nabla \varphi_h(\cdot  + k) \nonumber \\
&=&\int_{Q}  A_{\mathrm{per}} \left( \nabla w^{0,h}_p + p \right)\cdot \nabla \left( \sum_{|k|_{\infty} \leq N} \varphi_h(\cdot  + k) \right). \nonumber
\end{eqnarray}
Let $$\theta_h  := \sum_{|k|_{\infty} \leq N} \varphi_h(\cdot  + k). $$
Observe that $\theta_h$ is $Q$-periodic, and that its restriction to $Q$ is in the $\mathbb{P}_1$-Lagrange finite elements space built from $\mathcal{T}_h^{(Q)}$. In addition, we have 
$$\int_{Q_N} \theta_h = 0.$$
Thus $\theta_h \in V^h_{\mathrm{per}}(Q)$.
Hence by definition of $w^{0,h}_p$, we obtain, for any $\varphi_h \in V^h_\mathrm{per}(Q_N)$,
$$\int_{Q_N}  A_{\mathrm{per}} \left( \nabla w^{0,h}_p + p \right)\cdot \nabla \varphi_h = \int_{Q}  A_{\mathrm{per}}  \left(\nabla w^{0,h}_p+ p \right) \cdot \nabla \theta_h = 0.$$ In addition, we have $ w^{0,h}_p \in V^h_\mathrm{per}(Q) \subset V^h_\mathrm{per}(Q_N)$. As a consequence, $w^{0,h}_p$ is solution to~\eqref{eq:correcteur-0-hN-bis}. We conclude that $w_p^{0,h,N}(\cdot,\omega)=w^{0,h}_p$.

\paragraph{\underbar{Step 2:}} we introduce the function
\begin{equation}
\label{eq:def-v}
v^{\eta,h,N}_p = \frac{w^{\eta,h,N}_p - w^{0,h}_p}{\eta}. \nonumber
\end{equation}
We want to prove that $\nabla v^{\eta,h,N}_p $ is bounded in $\left( L^2(Q_N) \right)^d$ almost surely and independently from $\eta,h,N$ and $\omega$, and that it converges to the gradient of a random function, namely $w^{1,h,N}_p$ solution to~\eqref{eq:correcteur-1-hN}.
Substracting equation~(\ref{eq:correcteur-0-hN-bis}) from~(\ref{eq:diff-0-hN}), we obtain that $v^{\eta,h,N}_p(\cdot,\omega) \in V^h_\mathrm{per} (Q_N)$ is such that, for any $\varphi_h \in V_h^{\mathrm{per}} (Q_N),$
\begin{equation}
\label{eq:Cas1-PC-P0-N-v1}
\int_{Q_N} A_\mathrm{per} \nabla v_p^{\eta,h,N}(\cdot , \omega) \cdot \nabla \varphi_h + \int_{Q_N} \left( A_1(\cdot,\omega) + \eta^{-1} R_\eta(\cdot,\omega) \right) \left(p+ \nabla w_p^{\eta,h,N}(\cdot , \omega) \right)\cdot \nabla \varphi_h =0.
\end{equation}
Choosing $\varphi_h = v^{\eta,h,N}_p (\cdot,\omega)$ as test function in~\eqref{eq:Cas1-PC-P0-N-v1} and using~\eqref{eq:hyp-c-unif-0}, we have
\begin{eqnarray}
\label{eq:ineg-v}
\gamma	 \left\| \nabla v_p^{\eta,h,N}(\cdot , \omega)\right\|_{\left( L^2 (Q_N) \right)^d} &\leq & \left( \int_{Q_N} \left|\left(A_1(\cdot,\omega)  + \eta^{-1} R_\eta(\cdot,\omega) \right) p\right|^2  \right)^{\frac{1}{2}} \nonumber \\
&& \quad +  \left( \int_{Q_N} \left|\left(A_1(\cdot,\omega)  + \eta^{-1} R_\eta(\cdot,\omega) \right) \nabla w^{\eta,h,N}_p(\cdot,\omega) \right|^2  \right)^{\frac{1}{2}}, \nonumber \\
\label{eq:V}
&\leq&   \left( \left\|A_1 \right\|_{\left(L^\infty(Q \times \Omega)\right)^{d \times d}}  +  \left\| \eta^{-1} R_\eta (\cdot,\omega)\right\|_{\left(L^\infty(Q_N)\right)^{d \times d}} \right)\nonumber \\
&& \quad \times \left( \left\|\nabla w^{\eta,h,N}_p (\cdot,\omega) \right\|_{\left(L^2(Q_N)\right)^d} + |Q_N|^{\frac{1}{2}} \right).
\end{eqnarray}
Using~\eqref{eq:dl-H3}, we have, when $|\eta| \leq 1$, 
\begin{equation}
\label{eq:CR}
\left\| \eta^{-1} R_\eta(\cdot,\omega) \right\|_{\left(L^\infty(Q_N)\right)^{d \times d}} \leq \left\| \eta^{-1} R_\eta\right\|_{\left(L^\infty(Q \times \Omega) \right)^{d \times d}} \leq C_R |\eta|  \leq C_R.
\end{equation}
We then deduce from \eqref{eq:V}, \eqref{eq:Cw} and~\eqref{eq:CR} that, for any $|\eta| \leq 1$,
\begin{eqnarray}
\label{eq:Cv}
\left\| \nabla v^{\eta,h,N}_p (\cdot , \omega) \right\|_{\left(L^2 (Q_N)\right)^d} &\leq & \left( \left\|A_1 \right\|_{\left(L^\infty(\Omega \times Q)\right)^{d \times d}}  +  C_R \right) \left( \frac{C_w}{\gamma} + \frac{1	}{\gamma}   \right)|Q_N|^{\frac{1}{2}} \nonumber \\
&\leq& C_v |Q_N|^{\frac{1}{2}},
\end{eqnarray}
where 
\begin{equation}
\label{eq:def-Cv}
C_v = \left( \left\|A_1 \right\|_{\left(L^\infty(\Omega \times Q)\right)^{d \times d}}  +  C_R \right) \left( \frac{M}{\gamma^2} + \frac{1}{\gamma} \right) |Q_N|^{\frac{1}{2}} \nonumber
\end{equation}
is a constant independent of $\eta$, $h$, $N$ and $\omega$. Thus
$ \nabla v^{\eta,h,N}_p (\cdot,\omega)$ is also bounded independently of $\eta$ in $\left(L^2(Q_N)\right)^d$ almost surely.\\

As above, we deduce that there almost surely exists $w^{1,h,N}_p (\cdot,\omega) \in V^h_\mathrm{per}(Q_N)$ such that, up to extracting a subsequence,
\begin{equation}
\nabla v^{\eta,h,N}_p (\cdot,\omega) \mathop{\longrightarrow}_{\eta \rightarrow 0}  \nabla w^{1,h,N}_p (\cdot,\omega) \ \mbox{in} \ \left(L^2(Q_N)\right)^d.\nonumber
\end{equation}
Using the same arguments as in Step 1, we can pass to the limit $\eta \rightarrow 0$ in~\eqref{eq:Cas1-PC-P0-N-v1}. We obtain that $w^{1,h,N}_p(\cdot,\omega)$ satisfies
\begin{equation}
\label{eq:Cas1-P0-P1-N-FV}
\forall \varphi \in V_h^{\mathrm{per}} (Q_N), \ \ \int_{Q_N} A_{\mathrm{per}} \nabla w^{1,h,N}_p (\cdot,\omega) \cdot \nabla \varphi_h + \int_{Q_N} A_{1}(\cdot,\omega)\left(\nabla w_p^{0,h} +p \right) \cdot \nabla \varphi_h =0, \nonumber
\end{equation}
almost surely. We thus recover~\eqref{eq:correcteur-1-hN}.

\paragraph{\underbar{Step 3:}} Our purpose is now to obtain a second order approximation of $w_p^{\eta,h,N}$. To this end, we define 
\begin{equation}
\label{eq:def-t}
z^{\eta,h,N}_p = \frac{w_p^{\eta,h,N} - w^{0,h}_p - \eta w^{1,h,N}_p}{\eta^2} \in V^h_\mathrm{per} (Q_N). \nonumber
\end{equation}
Our goal is to prove that $\nabla z^{\eta,h,N}_p$ is bounded in $\left(L^2(Q_N)\right)^d$ almost surely and independently from $\eta,N,h$ and $\omega$. To do so, we will need
the previous estimates~\eqref{eq:Cw} and~\eqref{eq:Cv} from Steps 1 and 2 respectively.
Using~(\ref{eq:diff-0-hN}), (\ref{eq:correcteur-0-hN-bis}) and~(\ref{eq:correcteur-1-hN}) we see that, for any $\varphi_h \in V_h^{\mathrm{per}} (Q_N),$
\begin{equation}
\label{eq:Cas1-PC-P0-P1-N-v1}
\int_{Q_N} A_\mathrm{per} \nabla z_p^{\eta,h,N}(\cdot , \omega) \cdot \nabla \varphi_h + \int_{Q_N} A_1(\cdot,\omega) \nabla v_p^{\eta,h,N}(\cdot , \omega)\cdot \nabla \varphi_h 
+\int_{Q_N}\eta^{-2} R_\eta(\cdot,\omega)  \left(p+ \nabla w_p^{\eta,h,N}(\cdot , \omega) \right)\cdot \nabla \varphi_h =0 .
\end{equation}
We choose $\varphi_h = z^{\eta,h,N}_p(\cdot ,\omega)$ as test function in~\eqref{eq:Cas1-PC-P0-P1-N-v1}. Using~\eqref{eq:hyp-c-unif-0}, we obtain
\begin{eqnarray}
\label{eq:ineg-t}
\gamma \left\| \nabla z^{\eta,h,N}_p (\cdot,\omega) \right\|_{\left(L^2 (Q_N)\right)^d} &\leq& \left( \int_{Q_N} \left|A_1 \nabla v_p^{\eta,h,N}(\cdot , \omega) \right|^2 \right)^{\frac{1}{2}} + 
\left( \int_{Q_N} \left|\eta^{-2} R_\eta  \left(p+ \nabla w_p^{\eta,h,N}(\cdot , \omega) \right)\right|^2 \right)^{\frac{1}{2}}, \nonumber\\
&\leq&  \left\|A_1 \right\|_{\left(L^\infty(Q \times \Omega)\right)^{d \times d}} \left\|\nabla v_p^{\eta,h,N}(\cdot , \omega) \right\|_{\left(L^2(Q_N)\right)^d} \nonumber  \\
&& \quad + \left\|\eta^{-2} R_\eta(\cdot,\omega) \right\|_{\left(L^\infty(Q_N)\right)^{d \times d}} \left( \left\|\nabla w_p^{\eta,h,N}(\cdot , \omega) \right\|_{\left(L^2(Q_N)\right)^d} + \left| Q_N\right|^{\frac{1}{2}}  \right).
\end{eqnarray}
In view of~\eqref{eq:dl-H3} we have that, for any $\eta \leq 1$,
 $$ \left\| \eta^{-2} R_\eta (\cdot,\omega)\right\|_{\left(L^\infty(Q_N)\right)^{d \times d}} \leq \left\| \eta^{-2} R_\eta\right\|_{\left(L^\infty(Q \times \Omega) \right)^{d \times d}}\leq  C_R.$$
We deduce from~\eqref{eq:ineg-t}, \eqref{eq:Cw} and \eqref{eq:Cv} that
\begin{eqnarray}
\left\| \nabla z^{\eta,h,N}_p (\cdot,\omega)\right\|_{\left(L^2 (Q_N)\right)^d} &\leq&  \frac{\left\|A_1 \right\|_{\left(L^\infty(\Omega \times Q)\right)^{d \times d}}}{\gamma} C_v |Q_N|^{\frac{1}{2}} + \frac{C_R}{\gamma} (C_w + 1) \left|Q_N \right|^{\frac{1}{2}} \nonumber \\
\label{eq:Ct}
&\leq& C_z |Q_N|^{\frac{1}{2}},
\end{eqnarray}
with (see~\eqref{eq:Cw} and~\eqref{eq:Cv})
\begin{equation}
%\label{eq:Ct}
C_z = \frac{\left\|A_1 \right\|_{\left(L^\infty(\Omega \times Q)\right)^{d \times d}}}{\gamma}  \left( \left\|A_1 \right\|_{\left(L^\infty(\Omega \times Q)\right)^{d \times d}}  +  C_R \right) \left( \frac{M}{\gamma^2} + \frac{1}{\gamma} \right)
+ \frac{C_R}{\gamma} \left( \frac{M}{\gamma} +  1 \right). \nonumber
\end{equation}
We observe that $C_z$ is independent of $N,h,\eta$ and $\omega$. This concludes the proof of the first assertion of the proposition, namely the bound~\eqref{eq:cas1-dl-w-hN}.

\paragraph{\underbar{Step 4:}} We now prove the second assertion in the statement of Proposition 2.1, namely the expansion of the approximated homogenized matrix. Using~\eqref{eq:homog-matrix-Nh} and~\eqref{eq:dl}, we have for $1 \leq i,j\leq d$
\begin{eqnarray}
\left[ A^{\star,h,N}_{\eta}\right]_{ij}(\omega) &=& \frac{1}{|Q_N|} \int_{Q_N} e_i^T A_\mathrm{per} \left( \nabla w^{\eta,h,N}_{e_j}(\cdot,\omega) + e_j \right)  \eta \frac{1}{|Q_N|} \int_{Q_N} e_i^T A_1 (\cdot,\omega) \left( \nabla w^{\eta,h,N}_{e_j}(\cdot,\omega) + e_j \right) \nonumber \\
&+& \frac{1}{|Q_N|} \int_{Q_N} e_i^T R_\eta (\cdot,\omega) \left( \nabla w^{\eta,h,N}_{e_j}(\cdot,\omega) + e_j \right). \nonumber
\end{eqnarray}
Using definitions~\eqref{eq:a-0-hN},~\eqref{eq:a-1-hN},~\eqref{eq:def-v} and~\eqref{eq:def-t}, we write 
\begin{eqnarray}
\label{eq:A-A0-A1}
\eta^{-2} \left( \left[ A^{\star,h,N}_{\eta}\right]_{ij}(\omega) - \left[ A^{\star,h}_\mathrm{per}\right]_{ij} - \eta \left[ A_1^{\star,h,N}\right]_{ij}(\omega) \right) &=&\frac{1}{|Q_N|} \int_{Q_N} e_i^T A_\mathrm{per} \nabla z^{\eta,h,N}_{e_j}(\cdot,\omega) \nonumber \\
&+&  \frac{1}{|Q_N|} \int_{Q_N} e_i^T A_1 (\cdot,\omega) \nabla v^{\eta,h,N}_{e_j}(\cdot,\omega)  \nonumber \\ 
&+& \frac{1}{|Q_N|} \int_{Q_N} e_i^T \eta^{-2} R_\eta (\cdot,\omega)  \left( \nabla w^{\eta,h,N}_{e_j}(\cdot,\omega) + e_j \right). \nonumber
\end{eqnarray}
Applying Cauchy Schwarz inequality, and the bound~\eqref{eq:dl-H3}, we obtain that, for any $|\eta| \leq 1$,
\begin{eqnarray}
\label{eq:borne-a-error-2}
\eta^{-2} \left|\left[ A^{\star,h,N}_{\eta}\right]_{ij}(\omega) - \left[ A^{\star,h}_\mathrm{per}\right]_{ij} - \eta \left[ A_1^{\star,h,N} \right]_{ij} (\omega)\right| & \leq & \frac{1}{|Q_N|}  \left\|A_\mathrm{per} \right\|_{\left(L^\infty(Q)\right)^{d \times d}} |Q_N|^{\frac{1}{2}} \left\| \nabla z^{\eta,h,N}_{e_j} (\cdot,\omega) \right\|_{\left(L^2 (Q_N)\right)^d} \nonumber \\
 & + & \frac{1}{|Q_N|} \left\|A_1 \right\|_{\left(L^\infty(Q \times \Omega)\right)^{d \times d}} |Q_N|^{\frac{1}{2}} \left\| \nabla v^{\eta,h,N}_{e_j} (\cdot,\omega) \right\|_{\left(L^2 (Q_N)\right)^d} \nonumber \\
&+&\frac{C_R}{|Q_N|} \left(  |Q_N|^{\frac{1}{2}} \left\| \nabla w^{\eta,h,N}_{e_j} (\cdot,\omega) \right\|_{\left(L^2 (Q_N)\right)^d} + |Q_N| \right). \nonumber
\end{eqnarray}
Thus using~\eqref{eq:Cw},~\eqref{eq:Cv} and~\eqref{eq:def-Ct}, we deduce that
\begin{equation}
\label{eq:borne-a-error-2-bis}
\eta^{-2} \left|\left[ A^{\star,h,N}_{\eta}\right]_{ij}(\omega) - \left[ A^{\star,h}_\mathrm{per}\right]_{ij} - \eta \left[ A_1^{\star,h,N} \right]_{ij}(\omega) \right|
\leq M C_z + \left\|A_1 \right\|_{\left(L^\infty(\Omega \times Q)\right)^{d \times d}} C_v + C_R \left(1 + C_w \right). \nonumber
\end{equation}
Recall that $C_w$, $C_v$ and $C_z$ are all independent from $\eta$, $N$, $h$ and $\omega$. The estimate~\eqref{eq:cas1-dl-a-hN} is thus proved. This concludes the proof of Proposition~\ref{prop:cas1-dl-hN}.
\end{preuve21}

\begin{remarque}

The same proof applies under weaker hypotheses. Indeed, suppose that the expansion $A_\eta=A_{\mathrm{per}} + \eta A_1 + O(\eta^2)$ holds in $\left(L^\infty (\RR^d) \right)^{d \times d}$, almost surely. This means that~\eqref{eq:dl-H3} is replaced by
\begin{equation}
\label{eq:dl-H3-ps}
 \left\| \eta^{-2} R_\eta (\cdot,\omega)\right\|_{\left( L^\infty (\RR^d)\right)^{d \times d} } \leq C_R(\omega) \ \mbox{almost surely}, \nonumber
\end{equation}
where $C_R$ is now a random variable. We suppose in addition that the random variable $C_R$ satisfies $\EE \left(|C_R|^q \right)^{\frac{1}{q}} \leq \bar{C}_R$.
Under this weaker assumption, we have the following estimates (compare with ~\eqref{eq:cas1-dl-w-hN} and~\eqref{eq:cas1-dl-a-hN}):
\begin{eqnarray}
\label{eq:cas1-dl-w-hN-ps}
 \eta^{-2}\EE \left(  \left\| 
\nabla w^{\eta,h,N}_p  - \nabla w^{0,h}_p - \eta \nabla
w^{1,h,N}_p  \right\|_{\left(L^2(Q_N)\right)^{d}}^q \right)^{\frac{1}{q}} &\leq& \bar{C}_z \sqrt{| Q_N |}, \nonumber \\
\label{eq:dl-a-hN-ps}
 \eta^{-2} 
\EE \left( \left| A^{\star,h,N}_{\eta} - A_\mathrm{per}^{\star,h} - \eta A^{\star,h,N}_1  \right|^q \right)^{\frac{1}{q}}
&\leq& \bar{C}_A, \nonumber
\end{eqnarray}
with
\begin{eqnarray}
\label{eq:Ct-ps}
\bar{C}_z&=& \EE \left( \left| C_z \right|^q \right)^{\frac{1}{q}}, \nonumber \\
C_z(\omega)&=&\frac{\left\|A_1 \right\|_{\left(L^\infty(\Omega \times Q)\right)^{d \times d}}}{\gamma}  \left( \left\|A_1 \right\|_{\left(L^\infty(\Omega \times Q)\right)^{d \times d}}  +  C_R(\omega) \right) \left( \frac{M}{\gamma^2} + \frac{1}{\gamma} \right)
+ \frac{C_R(\omega)}{\gamma} \left( \frac{M}{\gamma} +  1 \right), \nonumber \\
\label{eq:CA-ps}
\bar{C}_A &=& \EE \left( \left| M C_z + \left\|A_1 \right\|_{\left(L^\infty(\Omega \times Q)\right)^{d \times d}} C_v + C_R \left(1 + C_w \right) \right|^q \right)^{\frac{1}{q}},\nonumber \\
C_v (\omega) &=& \left( \left\|A_1 \right\|_{\left(L^\infty(\Omega \times Q)\right)^{d \times d}}  +  C_R(\omega) \right) \left( \frac{M}{\gamma^2} + \frac{1}{\gamma} \right) |Q_N|^{\frac{1}{2}}. \nonumber
\end{eqnarray}
\end{remarque}

\begin{remarque}
In Proposition 2.1, the key hypothesis is~\eqref{eq:dl-H3}. This is the assumption from which we obtain inequalities~\eqref{eq:Cv} and~\eqref{eq:Ct} (from~\eqref{eq:ineg-v} and~\eqref{eq:ineg-t} respectively). Let us focus on~\eqref{eq:ineg-t}. Looking at the first line, the term we have to control reads:
$$
\left( \int_{Q_N} \left|\eta^{-2} R_\eta   \nabla w_p^{\eta,h,N}(\cdot , \omega)\right|^2 \right)^{\frac{1}{2}}.
$$
Using~\eqref{eq:dl-H3}, we have
$$
\left( \int_{Q_N} \left|\eta^{-2} R_\eta   \nabla w_p^{\eta,h,N}(\cdot , \omega)\right|^2 \right)^{\frac{1}{2}} \leq C_R \left( \int_{Q_N} \left| \nabla w_p^{\eta,h,N}(\cdot , \omega)\right|^2 \right)^{\frac{1}{2}},
$$ 
and using the bounds on $\nabla w_p^{\eta,h,N}$ previously obtained, we can conclude that $\nabla z^{\eta,h,N}_p$ is bounded in $\left(L^\infty \left(\Omega;L^2(Q_N) \right)  \right)^{d \times d}$. Suppose now that we work within a different framework, say the original expansion of $A_\eta$ holds in $\left(L^\infty \left(\RR^d;L^2(\Omega) \right) \right)^{d \times d}$. Hypothesis~\eqref{eq:dl-H3} is replaced by $$ \underset{x \in \reels^d}{\mathrm{Ess}\sup} \left( \EE \left( \left|R_\eta(x,\cdot) \right|^2 \right)^{\frac{1}{2}} \right) \leq C_R.$$ As a consequence, we cannot expect $\nabla z^{\eta,h,N}_p$ to be bounded in $\left(L^\infty \left(\Omega;L^2(Q_N) \right)  \right)^{d \times d}$. However, there may be a means to recover boundedness, but in a different space, namely $\left(L^2(Q_N \times \Omega)\right)^{d \times d}$. To do so, one needs $\nabla w^{\eta,h,N}_p \in \left(L^\infty \left(\Omega;L^\infty(Q_N) \right)  \right)^{d \times d}$, which requires stronger regularity hypotheses on~$A_\eta$.
\end{remarque}

\subsection{Convergence with respect to $h$ and $N$} 

In Proposition~\ref{prop:cas1-dl-hN}, we have obtained bounds for quantities defined at the discrete level, namely after truncation and finite elements discretization. We now study the limit of~\eqref{eq:cas1-dl-w-hN} and~\eqref{eq:cas1-dl-a-hN} when $h \rightarrow 0$, and next $N \rightarrow +\infty$.

\subsubsection{Convergence as $h \rightarrow 0$}

First, let us define $w_p^{\eta,N}$, $w^0_p $ and $w_p^{1,N}$ solutions to the problems
\begin{eqnarray} 
\label{eq:correcteur-N}
&&\left\{
\begin{array}{l}
\mbox{Find} \ 
w_p^{\eta,N}(\cdot,\omega) \in H^1_{\rm per}(Q_N) \ \ 
\mbox{such that,} \ 
\\
\displaystyle{\forall \varphi \in H^1_{\rm per}(Q_N), \ \ \int_{Q_N}
  A_\eta(\cdot,\omega)
\left( p +  \nabla w_p^{\eta,N}(\cdot,\omega)\right) \cdot \nabla \varphi = 0 \ \ \mbox{almost surely,}
}
\end{array}
\right.
\\
\label{eq:correcteur-0}
& &\left\{
\begin{array}{l}
\mbox{Find} \ 
w_p^{0} \in H^1_{\rm per}(Q) \ \ 
\mbox{such that,} \ 
\\
\displaystyle{\forall \varphi \in H^1_{\rm per}(Q), \ \ \int_{Q}
  A_{\mathrm{per}}
\left( p +  \nabla w_p^{0}\right) \cdot \nabla \varphi = 0,
}
\end{array}
\right.
\\
\label{eq:correcteur-1-N}
& &\left\{
\begin{array}{l}
\mbox{Find} \ 
w_p^{1,N}(\cdot,\omega) \in H^1_{\rm per}(Q_N) \ \ 
\mbox{such that,} \ 
\\
\displaystyle{\forall \varphi \in H^1_{\rm per}(Q_N), \ \ \int_{Q_N} A_{\mathrm{per}} \nabla w_p^{1,N}(\cdot,\omega) \cdot \nabla \varphi + \int_{Q_N}
  A_1(\cdot,\omega)
\left( p +  \nabla w_p^{0}\right) \cdot \nabla \varphi =0 \ \ \mbox{almost surely,}
}
\end{array}
\right.
\end{eqnarray}
respectively, where $H^1_\mathrm{per}(Q)$ denotes the closure of $C^\infty_\mathrm{per}(Q)$, the space of infinitly derivable periodic functions, with respect to the $H^1$-norm. We also define the matrices $A^{\star,N}_\eta$, $A_\mathrm{per}^{\star}$ and $A^{\star,N}_1$ by
\begin{eqnarray}
\label{eq:a-N}
\forall 1 \leq i,j \leq d, & & \left[A^{\star,N}_\eta \right]_{ij} (\omega) =\frac{1}{|Q_N|} \int_{Q_N} e_i^T A_\eta (\cdot,\omega) \left(e_j + \nabla w_{e_j}^{\eta,N}(\cdot,\omega) \right),\\
\label{eq:dl-a0-N}
&  & \left[ A_\mathrm{per}^{\star}\right] _{ij} = \int_Q e_i^T A_\mathrm{per} \,
\left( e_j + \nabla w_{e_j}^{0} \right), \\
\label{eq:dl-a1-N}
& &\left[ A^{\star,N}_1\right] _{ij}  (\omega) = \frac{1}{|Q_N|} \int_{Q_N} e_i^T A_\mathrm{per} \nabla w^{1,N}_{e_j} (\cdot,\omega)
+ \frac{1}{|Q_N|} \int_{Q_N} e_i^T A_1 (\cdot,\omega) \left( \nabla w^{0}_{e_j} + e_j \right). 
\end{eqnarray}
We now prove the following result, which is a direct consequence of Proposition~\ref{prop:cas1-dl-hN}.
\begin{corollaire}
\label{prop:cas1-dl-N}
Suppose that $A_\eta$  is a symmetric matrix that satisfies (\ref{eq:hyp-bc-unif}) and is stationary in the sense of~\eqref{eq:stationnarite-disc}. Suppose, in addition, that it satisfies~(\ref{eq:dl}), (\ref{eq:dl-H1}), (\ref{eq:dl-H2}). We assume that (\ref{eq:dl-H3}) holds, namely the second order error is $O(\eta^2)$ in $\left(L^\infty (Q \times \Omega) \right)^{d \times d}$. Then there exists a constant $C$, independent of $N$, $\eta$ and $\omega$ such that, for $|\eta| \leq 1$,
\begin{equation}
\label{eq:cas1-dl_w_N}
 \eta^{-2} \left\| 
\nabla w^{\eta,N}_p(\cdot,\omega) - \nabla w^{0}_p - \eta \nabla
w^{1,N}_p(\cdot,\omega) \right\|_{\left(L^2(Q_N)\right)^{d}} \leq C \sqrt{| Q_N |} \ \ \mbox{almost surely},
\end{equation}
where $w^{\eta,N}_p$, $w^{0}_p$ and $w^{1,N}_p$ are solutions
to~(\ref{eq:correcteur-N}), (\ref{eq:correcteur-0})
and~(\ref{eq:correcteur-1-N}), respectively, and such that
\begin{equation}
\label{eq:dl-a-N}
 \eta^{-2} 
\left| A^{\star,N}_{\eta}(\omega) - A_\mathrm{per}^{\star} - \eta A^{\star,N}_1(\omega) \right|
\leq C \ \ \mbox{almost surely},
\end{equation}
where $A^{\star,N}_\eta$, $A_\mathrm{per}^{\star}$ and $A^{\star,N}_1$ are defined by~\eqref{eq:a-N},~\eqref{eq:dl-a0-N} and ~\eqref{eq:dl-a1-N} respectively. 
\end{corollaire}

\begin{preuve} Using~\eqref{eq:cas1-dl-w-hN}, remark that
\begin{eqnarray}
\label{eq:ineq}
 \eta^{-2} \left\| 
\nabla w^{\eta,N}_p(\cdot,\omega) - \nabla w^{0}_p - \eta \nabla
w^{1,N}_p(\cdot,\omega) \right\|_{\left(L^2(Q_N)\right)^{d}} &\leq&
\eta^{-2} \left\|\nabla w^{\eta,h,N}_p(\cdot,\omega) - \nabla w^{\eta,N}_p(\cdot,\omega)\right\|_{\left(L^2(Q_N)\right)^{d}} \nonumber \\
&+& \eta^{-2} \left\| \nabla w^{0,h}_p - \nabla w^0_p \right\|_{\left(L^2(Q_N)\right)^{d}} \nonumber \\
&+& \eta^{-1} \left\| \nabla
w^{1,h,N}_p(\cdot,\omega)  - \nabla w^{1,N}_p (\cdot,\omega) \right\|_{\left(L^2(Q_N)\right)^{d}} \\
&+& C \sqrt{| Q_N |} \nonumber,
\end{eqnarray}
where $C$ is independent of $h$, $N$, $\omega$ and $\eta$.
Using standard properties of finite element approximations, we have that
\begin{itemize}
\item $\underset{h \rightarrow 0}{\lim} \left\|\nabla w^{\eta,h,N}_p(\cdot,\omega) - \nabla w_p^{\eta,N}(\cdot,\omega) \right\|_{\left(L^2(Q_N)\right)^d} = 0$ almost surely,
\item $\underset{h \rightarrow 0}{\lim} \left\|\nabla w^{1,h,N}_p(\cdot,\omega) - \nabla w_p^{1,N}(\cdot,\omega) \right\|_{\left(L^2(Q_N)\right)^d} = 0$ almost surely,
\item $\underset{h \rightarrow 0}{\lim} \left\|\nabla w^{0,h}_p - \nabla w_p^{0} \right\|_{\left(L^2(Q_N)\right)^d} = 0$.
\end{itemize}
Passing to the limit $h \rightarrow 0$ in~\eqref{eq:ineq}, we obtain~\eqref{eq:cas1-dl_w_N}. The proof of~\eqref{eq:dl-a-N} follows the same line.
% Making $h$ tend to $0$, we obtain that the expansion holds in $\left(L^2(Q_N) \right)^d$ almost surely. But once again, since the constant $C$ is independent
% from $\omega$ (basically it is the same as in Proposition~\ref{prop:cas1-dl-hN}),  is proved.
\end{preuve}

\subsubsection{Convergence as $N \rightarrow \infty$}

We now study the limit of~\eqref{eq:dl-a-N} as $N \rightarrow + \infty$. From~\cite{BourgeatPiatnitski} we already know that 
\begin{equation}
\label{eq:conv-a-N}
\underset{N \rightarrow + \infty}{\lim} A^{\star,N}_\eta (\omega) = A^\star_\eta \ \mbox{almost surely}, 
\end{equation}
where the exact homogenized matrix $A^\star_\eta$ is defined by~\eqref{eq:Cas1-Astar}. We now turn to $A^{\star,N}_1$ whose limit when $N \rightarrow + \infty $ is given by the following lemma.
\begin{lemme} 
Suppose that the matrices $A_\eta$, $A_1$ and $A_\mathrm{per}$ are symmetric, then the matrix $A^{\star,N}_1 (\omega)$ defined by~\eqref{eq:dl-a1-N} satisfies
\begin{equation}
\label{eq:conv-a1-N}
\underset{N \rightarrow + \infty}{\lim} A^{\star,N}_1 (\omega) = A_1^\star\ \mbox{almost surely}, 
\end{equation}
where the deterministic matrix $A^\star_1$ is given by
\begin{equation}
\left[ A^{\star}_1\right] _{ij}  (\omega) = \int_{Q} \left(e_i + \nabla w^0_{e_i} \right)^T  \esp \left( A_1 \right) \left( e_j + \nabla w^{0}_{e_j} \right), \nonumber
\end{equation} 
where $w^{0}_p$ is the unique solution (up to an additive constant) to~\eqref{eq:correcteur-0}.
\end{lemme}

\begin{preuve}
We observe that
\begin{eqnarray}
\left[ A^{\star,N}_1\right] _{ij}  (\omega) &=& \frac{1}{|Q_N|} \int_{Q_N} e_i^T A_\mathrm{per} \nabla w^{1,N}_{e_j} (\cdot,\omega)
+ \frac{1}{|Q_N|} \int_{Q_N} e_i^T A_1 (\cdot,\omega) \left( \nabla w^{0}_{e_j} + e_j \right) \nonumber \\
&=& \frac{1}{|Q_N|} \int_{Q_N} \left(e_i + \nabla w^0_{e_i} \right)^T \left(  A_\mathrm{per} \nabla w^{1,N}_{e_j} (\cdot,\omega) + A_1 (\cdot,\omega) \left( \nabla w^{0}_{e_j} + e_j \right) \right)
\end{eqnarray}
because $ w^0_{e_i} $ is $Q_N$-periodic and can therefore be chosen as a test function in~\eqref{eq:correcteur-1-N}. Since $A_\mathrm{per}$ is symmetric, we have,  using~\eqref{eq:correcteur-0}, that
$$\int_{Q_N} \left(e_i + \nabla w^0_{e_i} \right)^T A_\mathrm{per} \nabla w^{1,N}_{e_j} (\cdot,\omega) = \int_{Q_N}  \left(\nabla w^{1,N}_{e_j}(\cdot,\omega)\right)^T A_\mathrm{per}  \left(e_i + \nabla w^0_{e_i} \right) = 0 .$$
Thus
$$ \left[ A^{\star,N}_1\right] _{ij}  (\omega)= \frac{1}{|Q_N|} \int_{Q_N} \left(e_i + \nabla w^0_{e_i} \right)^T A_1 (\cdot,\omega) \left( \nabla w^{0}_{e_j} + e_j \right),$$
and we conclude using the ergodic Theorem~\ref{theo:ergodic}, as $A_1$ is stationary and $\nabla w^0_p$ is periodic.
\end{preuve}
Using~\eqref{eq:conv-a-N} and~\eqref{eq:conv-a1-N}, we can pass to the limit $N \rightarrow +\infty$ in~\eqref{eq:dl-a-N} and we obtain that
\begin{equation}
\label{eq:dl-a}
 \eta^{-2} 
\left| A^{\star}_{\eta} - A_\mathrm{per}^{\star} - \eta A^{\star}_1 \right|
\leq C,
\end{equation}
where $C$ is independent of $\eta$ and $\omega$ (note that $A^{\star}_{\eta}$, $ A_\mathrm{per}^{\star}$ and $ A^{\star}_1 $ are all deterministic matrices).  We thus recover the expansion of the exact \emph{deterministic} homogenized matrix $A^\star_\eta$ as given in \cite{MPRF}.

\section{Expansion of stochastic diffeomorphism (Model 2)}

We now focus on Model 2. The goal of this section is the same as that of Section 2. We  prove that the random second order error in the expansions of the gradient of corrector $\nabla w^{\eta,h,N}_p$ and the homogenized matrix $A^{\star,h,N}_\eta$ is bounded independently of $h$, $N$ and $\eta$ in some appropriate $L^p(\Omega)$ space. The functional space in which the expansion of the original diffeomorphism $\Phi_\eta$ holds is simple and somehow corresponds to the one considered for Model 1 in Section 2. 

\subsection{Hypotheses}

In this section, we consider the Model 2 mentioned above, where the random matrix $A_\eta$ in~\eqref{eq:edp-epsilon} writes
\begin{equation}
A_\eta (x,\omega) = A_\mathrm{per} \left(\Phi_\eta^{-1}(x,\omega) \right),
\end{equation}
where $\Phi_\eta$ is a \emph{stochastic diffeomorphism}, that satisfies conditions~\eqref{eq:hyp-diffeo-as},~\eqref{eq:hyp-diffeo-1},~\eqref{eq:hyp-diffeo-2} and~\eqref{eq:hyp-diffeo-3}. 
The periodic matrix $A_\mathrm{per}$ is supposed uniformly bounded and coercive:
\begin{equation}
\label{eq:hyp-bc-unif-per}
\left\{
\begin{array}{l}
\exists \gamma >0, \ \forall \xi \in \RR^d, \ \xi^T A_\mathrm{per}(x) \xi \geq \gamma |\xi|^2 \ \mbox{almost everywhere}, \\
\exists M >0 \mbox{ such that } \left\|A_\mathrm{per}  \right\|_{\left(L^\infty(Q)\right)^{d \times d}} \leq M.
\end{array}
\right.
\end{equation}
This model introduced by Blanc, Le Bris and Lions in~\cite{BlancLeBrisLions1} is not a particular case of the standard homogenization setting. Following~\cite{BlancLeBrisLions2} we now consider a weakly stochastic case, where the diffeomorphism $\Phi_\eta$ is close to the identity. More precisely, we suppose in the sequel that the following expansion holds in $\left(C^1 \left(\RR^d, L^\infty(\Omega) \right) \right)^{d}$
\begin{equation}
\label{eq:dl-cas2}
\Phi_\eta(x,\omega) = x + \eta \Psi(x,\omega) + \Theta_\eta(x,\omega) \ \ \mbox{with} \ \ \Theta_\eta=O(\eta^2),
\end{equation}
where $\Psi$ satisfies~\eqref{eq:hyp-diffeo-2} and~\eqref{eq:hyp-diffeo-3}.
This means that
\begin{eqnarray}
\label{eq:dl-H1-cas2}
\underset{\eta \rightarrow 0}{\lim} \left\|\Phi_\eta - \mbox{Id} \right\|_{C^1 \left(\RR^d, L^\infty(\Omega) \right)^{d}}&=&0,  \\
\label{eq:dl-H2-cas2}
\underset{\eta \rightarrow 0}{\lim} \left\| \eta^{-1} \left(\Phi_\eta - \mbox{Id}\right) -\Psi \right\|_{C^1 \left(\RR^d, L^\infty(\Omega) \right)^{d}}&=&0,
\end{eqnarray}
where $\mbox{Id}$ denotes the identity mapping. In addition, $\Theta_\eta=O(\eta^2) $ means that there exists a \emph{deterministic} constant $ C_\Theta $ independent of $\eta$ such that, when $|\eta| \leq 1 $, 
\begin{equation}
\label{eq:dl-H3-cas2}
\left\| \eta^{-2} \Theta_\eta  \right\|_{C^1 \left(\RR^d, L^\infty(\Omega) \right)^{d}} \leq C_\Theta. 
\end{equation}
In~\cite{BlancLeBrisLions2}, it has been shown  that, under these hypotheses, the exact homogenized matrix~\eqref{eq:Cas2-Astar} possesses an expansion $A^\star_\eta= A^\star_\mathrm{per} + \eta A^\star_1 + O(\eta^2)$, where $A^\star_\mathrm{per}$ and $A^\star_1$ are deterministic matrices that only involve solutions to deterministic problems posed on a bounded domain, in contrast to $A^\star_\eta$ (see~\eqref{eq:Cas2-PC-Linear} and~\eqref{eq:Cas2-Astar}). Here, as above, we will prove that the approximated homogenized matrix, obtained after truncation and discretization, also possesses an expansion in powers of $\eta$ and that the error at second order is bounded independently of the parameters of the discretization procedure.

\subsection{Formal expansion}

As in Section 2, we first present formally our main result.
We assume the following formal expansion on the approximated corrector, solution to~\eqref{eq:correcteur-hN-cas2}:
\begin{equation}
\label{eq:dl-w-cas2}
w_p^{\eta,h,N} = w^{0,h,N}_p + \eta w^{1,h,N}_p + r_p^{\eta,h,N},
\end{equation}
where $\nabla r_p^{\eta,h,N} = O (\eta^2)$ in some appropriate space. Inserting this expansion in~\eqref{eq:correcteur-hN-cas2}, we obtain that the function $w_p^{0,h,N}$ is in fact equal to $w_p^{0,h}$ solution to~\eqref{eq:correcteur-0-hN}.
The function $w_{p}^{1,h,N}(\cdot,\omega) \in  V^h_\mathrm{per}(Q_N)$ is the unique function such that, for any $\varphi_h \in V^h_\mathrm{per}(Q_N) $, we have
\begin{eqnarray}
\label{eq:correcteur-1-hN-cas2}
& &\int_{Q_N}  A_{\mathrm{per}}
 \nabla w_p^{1,h,N}(\cdot,\omega) \cdot \nabla \varphi_h + \int_{Q_N}  A_{\mathrm{per}}
\nabla \Psi(\cdot,\omega) \nabla w_p^{0,h,N}(\cdot,\omega) \cdot \nabla \varphi_h  \nonumber \\
&+&\int_{Q_N} \left( \mathrm{div} \Psi (\cdot,\omega) \mathbb{I}_d - \nabla \Psi(\cdot,\omega)^T \right) A_{\mathrm{per}}
\left( p +  \nabla w_p^{0,h,N}(\cdot,\omega)\right) \cdot \nabla \varphi_h =0 \ \ \mbox{almost surely}.
\end{eqnarray}
In addition, substituting $A_\eta$ and $w_p^{\eta,h,N}$  by their expansions into~\eqref{eq:a-hN-cas2}, we formally obtain that 
\begin{eqnarray}
\label{eq:dl-a-hN-cas2-1}
A^{\star,h,N}_\eta (\omega) = A^{\star,h}_\mathrm{per} + \eta A_1^{\star,h,N}(\omega) + O(\eta^2),
\end{eqnarray}
where the term of order zero is defined by~\eqref{eq:correcteur-0-hN} and~\eqref{eq:a-0-hN}, as in the first model. The term of order one is defined by
\begin{eqnarray}
\label{eq:a-1-hN-cas2}
\left[A^{\star,h,N}_1 \right]_{ij}(\omega)&=&-\left[A^{\star,h}_\mathrm{per} \right]_{ij} \frac{1}{|Q_N|}\int_{Q_N}  \mathrm{div} \Psi(\cdot,\omega) + \frac{1}{|Q_N|}\int_{Q_N} \mathrm{div} \Psi (\cdot,\omega) \left(e_i + \nabla w^{0,h}_{e_i} \right)^T A_{\mathrm{per}} e_j \nonumber \\
&+& \frac{1}{|Q_N|}\int_{Q_N} \left(\nabla w_{e_i}^{1,h,N}(\cdot,\omega) - \nabla \Psi (\cdot,\omega) \nabla w^{0,h}_{e_i} \right)^T A_{\mathrm{per}} e_j.
\end{eqnarray}
As in Section 2, we now make precise and rigorously justify the expansions~\eqref{eq:dl-w-cas2} and~\eqref{eq:dl-a-hN-cas2-1}.

\subsection{Main result}

The main result of this section is the following proposition, which is analogous to Proposition~\ref{prop:cas1-dl-hN}.
\begin{prop}
\label{prop:cas2-dl-hN}
Suppose that $\Phi_\eta$ satisfies (\ref{eq:hyp-diffeo-as}),~(\ref{eq:hyp-diffeo-1}),~(\ref{eq:hyp-diffeo-2}) and~(\ref{eq:hyp-diffeo-3}), and that it is a perturbation of the identity in the sense of (\ref{eq:dl-cas2}), (\ref{eq:dl-H1-cas2}), (\ref{eq:dl-H2-cas2}) and (\ref{eq:dl-H3-cas2}). We suppose in addition that the symmetric periodic matrix $A_\mathrm{per}$ satisfies~\eqref{eq:hyp-bc-unif-per}. Then there exists a constant $C$ independent of $\eta$, $\omega$, $N$ and $h$, such that, for sufficiently small values of $\eta$,
\begin{equation}
\label{eq:dl-w-hN-cas2}
 \eta^{-2} \left\| 
\nabla w^{\eta,h,N}_p(\cdot,\omega) - \nabla w^{0,h}_p - \eta \nabla
w^{1,h,N}_p(\cdot,\omega) \right\|_{\left(L^2(Q_N)\right)^{d}} \leq C \sqrt{| Q_N |} \ \ \mbox{almost surely}
\end{equation}
where $w^{\eta,h,N}_p$, $w^{0,h}_p$ and $w^{1,h,N}_p$ are solutions
to~(\ref{eq:correcteur-hN-cas2}), (\ref{eq:correcteur-0-hN})
and~(\ref{eq:correcteur-1-hN-cas2}), respectively, and such that
\begin{equation}
\label{eq:dl-a-hN-cas2}
 \eta^{-2} 
\left| A^{\star,h,N}_{\eta}(\omega) - A_\mathrm{per}^{\star,h} - \eta A^{\star,h,N}_1(\omega) \right|
\leq C \ \ \mbox{almost surely},
\end{equation}
where $A^{\star,h,N}_\eta$, $A^{\star,h}_\mathrm{per}$ and $A^{\star,h,N}_1$ are given by~\eqref{eq:a-hN-cas2},~\eqref{eq:a-0-hN} and~\eqref{eq:a-1-hN-cas2} respectively.
\end{prop}

The proof, which is detailed below, follows the same lines as that of Proposition~\ref{prop:cas1-dl-hN}. It makes use of the following lemma.
\begin{lemme}
\label{lemme-diffeo}
Suppose that $\Phi_\eta$ satisfies (\ref{eq:dl-cas2}), (\ref{eq:dl-H1-cas2}), (\ref{eq:dl-H2-cas2}) and (\ref{eq:dl-H3-cas2}). Then,
there exist $$ \Gamma_\eta \in \left(C^0 \left( \RR^d, L^\infty(\Omega) \right)\right)^{d \times d} \ \mbox{and} \ \sigma_\eta \in C^0 \left(\RR^d, L^\infty(\Omega)  \right),$$ such that
\begin{eqnarray}
\label{eq:dl-gradphi}
\left(\nabla  \Phi_\eta \right)^{-1}(\cdot,\omega)  &=& \mathbb{I}_d - \eta \nabla \Psi(\cdot,\omega) + \Gamma_\eta (\cdot,\omega) \ \mbox{almost surely}, \\
\label{eq:dl-detphi}
\mathrm{det}(\nabla \Phi_\eta)(\cdot,\omega) &=& 1 + \eta \mathrm{div} \Psi (\cdot,\omega) + \sigma_\eta (\cdot,\omega) \ \mbox{almost surely},
\end{eqnarray}
with, when $\eta$ is sufficiently small,
\begin{equation}
\label{eq:dl-gamma-sigma}
\left\|\Gamma_\eta \right\|_{C^0 \left( \RR^d, L^\infty(\Omega) \right)^{d \times d}} \leq C_\Gamma \eta^2 \ \ \mbox{and} \ \ \left\|\sigma_\eta \right\|_{C^0 \left( \RR^d, L^\infty(\Omega)  \right)} \leq C_\sigma \eta^2,
\end{equation}
for two deterministic constants $C_\Gamma$ and $C_\sigma$ independent of $\eta$. In addition, there exists $\eta_0$ such that, for all $|\eta| \leq \eta_0$, we have
\begin{equation}
\label{eq:dl-valpphi}
\forall  \xi \in \RR^d, \  \forall x \in \RR^d, \ \ \frac{1}{2} |\xi|^2 \leq \xi^T \left( \nabla \Phi_\eta (x,\omega) \right)^{-T} \cdot \left(\nabla \Phi_\eta (x,\omega) \right)^{-1} \xi \leq \frac{3}{2} |\xi|^2 \ \  \ \ \mbox{almost surely}.
\end{equation}
\end{lemme}

\begin{preuvebis}
The proofs of~\eqref{eq:dl-gradphi},~\eqref{eq:dl-detphi} and \eqref{eq:dl-gamma-sigma} are straightforward. We now prove~\eqref{eq:dl-valpphi}. First, let us denote by $\Lambda$ the application that associates to any symmetric matrix $M \in \mathbb{S}_d(\RR)$ its ordered eigenvalues $\Lambda(M) = (\lambda_1,\cdots,\lambda_d)$. This application is continuous on $\mathbb{S}_d(\RR)$. As a consequence, there exists a value $\delta_0$ such that
\begin{equation}
\label{lemma2-eq1}
\forall \delta \leq \delta_0, \ \left| M - \mathbb{I}_d \right| \leq \delta \Rightarrow \left|\Lambda(M) -(1,\cdots,1)\right| \leq \frac{1}{2}. 
\end{equation}
We infer from~\eqref{eq:dl-gradphi} that the following expansion holds in $ \left(C^0 \left( \RR^d, L^\infty(\Omega) \right) \right)^{d \times d}$: $$\left(\nabla  \Phi_\eta \right)^{-T} \left(\nabla  \Phi_\eta \right)^{-1}  = \mathbb{I}_d - \eta \left( \left( \nabla \Psi\right)^T + \nabla \Psi \right) + O(\eta^2).$$
As a consequence, there exists a \emph{deterministic} value $\eta_0$ such that
\begin{equation}
\label{lemma2-eq2}
\forall \eta \leq \eta_0, \ \forall x \in \RR^d, \ \ \left|\left(\nabla  \Phi_\eta(x,\omega) \right)^{-T} \left(\nabla  \Phi_\eta(x,\omega) \right)^{-1} - \mathbb{I}_d\right| \leq \delta_0, \ \ \ \mbox{almost surely}. 
\end{equation}
Collecting~\eqref{lemma2-eq1} and~\eqref{lemma2-eq2}, we deduce that 
$$\forall |\eta| \leq \eta_0,\ \forall x \in \RR^d, \  \left|\Lambda \left( \left(\nabla  \Phi_\eta(x,\omega) \right)^{-T} \left(\nabla  \Phi_\eta(x,\omega) \right)^{-1} \right) -(1,\cdots,1)\right| \leq \frac{1}{2},  \ \ \mbox{almost surely}.$$
Hence all the eigenvalues of $\left(\nabla  \Phi_\eta \right)^{-T} \left(\nabla  \Phi_\eta \right)^{-1} $ belong to the interval $\left[1/2,3/2\right]$. It implies~\eqref{eq:dl-valpphi}.
\end{preuvebis}

\begin{preuveter}
As announced above, the structure of the proof is the same as that of Proposition 2.1. First, we rigorously justify~\eqref{eq:dl-w-cas2}. Doing so, we obtain bounds on $\nabla w^{\eta,h,N}_p$, and on the error terms $\nabla v^{\eta,h,N}_p$ and $\nabla z^{\eta,h,N}_p$ that are defined below (Steps 1,2 and 3 respectively). These bounds are independent of $h$, $N$ and $\eta$. We use them at Step 4 to control the second order error of the expansion of the approximated homogenized matrix $A^{\star,h,N}_\eta$.

\paragraph{\underbar{Step 1:}} Our goal is first to prove that $ \nabla w^{\eta,h,N}_p(\cdot,\omega)$ is bounded in $\left(L^2(Q_N)\right)^d$, almost surely. We thus show that $w^{0,h,N}_p$, the limit of $ w^{\eta,h,N}_p$, is in fact the unique solution to~\eqref{eq:correcteur-0-hN}. Choosing $\varphi_h = w^{\eta,h,N}_p(\cdot,\omega)$ as test function in~\eqref{eq:correcteur-hN-cas2} and using~\eqref{eq:hyp-bc-unif-per}, we obtain that 
\begin{eqnarray}
%\label{eq:W}
\gamma \left\| \sqrt{\mathrm{det}(\nabla \Phi_\eta)(\cdot,\omega)} \left(\nabla \Phi_\eta(\cdot,\omega) \right)^{-1} \nabla w^{\eta,h,N}_p(\cdot,\omega) \right\|_{\left(L^2 (Q_N)\right)^{d}} 
 \leq  \left(\int_{Q_N}\mathrm{det}(\nabla \Phi_\eta) (\cdot,\omega)  |   A_\mathrm{per}  p |^2 \right)^{\frac{1}{2}}. \nonumber 
\end{eqnarray}
Using Lemma~\ref{lemme-diffeo} and~\eqref{eq:hyp-diffeo-1} we obtain
\begin{equation}
\label{prop-2-step1-eq1}
\gamma \left\| \sqrt{\mathrm{det}(\nabla \Phi_\eta)(\cdot,\omega)} \left(\nabla \Phi_\eta(\cdot,\omega) \right)^{-1} \nabla w^{\eta,h,N}_p(\cdot,\omega) \right\|_{\left(L^2 (Q_N)\right)^{d}} \geq \gamma \sqrt{\frac{\nu}{2}}\left\| \nabla w^{\eta,h,N}_p(\cdot,\omega) \right\|_{\left(L^2 (Q_N)\right)^{d}}.
\end{equation}
From~\eqref{eq:hyp-diffeo-2}, we know that $\left| \mathrm{det}(\nabla \Phi_\eta)(\cdot,\omega) \right|\leq d! (M')^d $. We thus have (recalling that $|p|=1$)
\begin{equation}
\label{prop-2-step1-eq2}
\left(\int_{Q_N} \mathrm{det}(\nabla \Phi_\eta) (\cdot,\omega)  |   A_\mathrm{per}  p |^2 \right)^{\frac{1}{2}} \leq \left( d! (M')^d\right)^{\frac{1}{2}} M |Q_N|^{\frac{1}{2}}.
\end{equation}
Collecting~\eqref{prop-2-step1-eq1} and~\eqref{prop-2-step1-eq2}, we obtain that
\begin{equation}
\label{eq:Cw-diffeo}
\left\|\nabla w^{\eta,h,N}_p(\cdot,\omega) \right\|_{\left(L^2 (Q_N)\right)^{d}} \leq C_{w} |Q_N|^{\frac{1}{2}},
\end{equation}
where $$C_{w} = \sqrt{\frac{2}{\nu}} \frac{\left( d! (M')^d\right)^{\frac{1}{2}} M}{\gamma}$$ is independent of $\eta$, $h$, $N$ and $\omega$. Hence, $\nabla w^{\eta,h,N}_p$ is bounded in $L^2 \left(Q_N \right)^d$, independently of $\eta$. As in the proof of Proposition~\ref{prop:cas1-dl-hN} we deduce that there exists a
function $ w^{0,h,N}_p(\cdot,\omega) \in V^{\mathrm{per}}_h(Q_N)$ such that
\begin{equation}
\label{eq:conv-w0}
\nabla w^{\eta,h,N}_p (\cdot,\omega) \mathop{\longrightarrow}_{\eta \rightarrow 0} \nabla w^{0,h,N}_p \ \mbox{in} \left(L^2(Q_N)\right)^{d} \ \mbox{almost surely.}
\end{equation}
We now insert in equation~\eqref{eq:correcteur-hN-cas2} the expansions~\eqref{eq:dl-gradphi} and~\eqref{eq:dl-detphi}. For any $\varphi_h \in V^h_\mathrm{per} (Q_N)$,
\begin{eqnarray}
\label{eq:PC-dev-diffeo}
0&= & \int_{Q_N}
  \mathrm{det}(\nabla \Phi_\eta(\cdot,\omega))(\nabla \Phi_\eta(\cdot,\omega))^{-T} A_{\mathrm{per}}
\left( p + (\nabla \Phi_\eta(\cdot,\omega))^{-1} \nabla w_p^{\eta,h,N}(\cdot,\omega)\right) \cdot \nabla \varphi_h  \nonumber \\
\label{eq:order0}
 &=& E^\eta_0(\omega) + \eta E^\eta_1(\omega) + \eta^2 E^\eta_2(\omega) + E^\eta_3(\omega) + F_\eta(\omega),
\end{eqnarray}
with
\begin{eqnarray}
E^\eta_0( \omega) &=& \int_{Q_N}  A_{\mathrm{per}}
\left( p +  \nabla w_p^{\eta,h,N}(\cdot,\omega)\right) \cdot \nabla \varphi_h \nonumber \\
E^\eta_1(\omega) &=& \int_{Q_N} \mathrm{div} \Psi (\cdot,\omega) A_{\mathrm{per}}
\left( p +  \nabla w_p^{\eta,h,N}(\cdot,\omega)\right) \cdot \nabla \varphi_h 
-  \int_{Q_N}  \nabla \Psi(\cdot,\omega)^T   A_{\mathrm{per}}
\left( p +  \nabla w_p^{\eta,h,N}(\cdot,\omega)\right) \cdot \nabla \varphi_h \nonumber \\
&-& \int_{Q_N}  A_{\mathrm{per}}
\nabla \Psi(\cdot,\omega) \nabla w_p^{\eta,h,N}(\cdot,\omega) \cdot \nabla \varphi_h, \\
E^\eta_2(\omega) &=& \int_{Q_N} \nabla \Psi(\cdot,\omega)^T A_{\mathrm{per}} \nabla \Psi(\cdot,\omega) \nabla w_p^{\eta,h,N}(\cdot,\omega) \cdot \nabla \varphi_h 
- \int_{Q_N} \mathrm{div} \Psi (\cdot,\omega) \nabla \Psi(\cdot,\omega)^T A_{\mathrm{per}}  \left( \nabla w_p^{\eta,h,N}(\cdot,\omega)  + p \right) \cdot \nabla \varphi_h \nonumber \\
&-& \int_{Q_N} \mathrm{div} \Psi (\cdot,\omega)  A_{\mathrm{per}} \nabla \Psi(\cdot,\omega) \nabla w_p^{\eta,h,N}(\cdot,\omega) \cdot \nabla \varphi_h, \\
E^\eta_3(\omega) &=& \int_{Q_N} \sigma_\eta (\cdot,\omega) A_{\mathrm{per}}
\left( p +  \nabla w_p^{\eta,h,N}(\cdot,\omega)\right) \cdot \nabla \varphi_h 
+ \int_{Q_N} \Gamma_\eta (\cdot,\omega)^T  A_{\mathrm{per}}  \left( \nabla w_p^{\eta,h,N}(\cdot,\omega)  + p \right) \cdot \nabla \varphi_h \nonumber \\
&+& \int_{Q_N}  A_{\mathrm{per}} \Gamma_\eta(\cdot,\omega) \nabla w_p^{\eta,h,N}(\cdot,\omega) \cdot \nabla \varphi_h,
\end{eqnarray}
and  $F_\eta(\omega)$ is the remainder term. Formally, $E_0^\eta$ is of order $\eta^0$, $\eta E_1^\eta$ is of order $\eta$, $\eta^2 E_2^\eta$ is of order $\eta^2$ (because it includes products of two quantities of the order $\eta$), $E_3^\eta$ is of order $\eta^2$ as well (it includes products of a quantity of order $\eta^2$ with a quantity of order $\eta^0$) and $F_\eta$ is of order $\eta^3$. We now prove rigorously boundedness and convergence results for each of the terms: $E^\eta_0, E^\eta_1,E^\eta_2,E^\eta_3$ and $F_\eta$. These estimates are very useful not only at Step 1 but also at Steps 2 and 3. \\

\paragraph{Convergences of $E^\eta_0$, $E^\eta_1$ and $E^\eta_2$, as $\eta \rightarrow 0$.}
We deduce from~\eqref{eq:conv-w0} that 
\begin{eqnarray}
\label{eq:conv-E0}
E^\eta_0(\omega) =\int_{Q_N}  A_{\mathrm{per}} \left( p +  \nabla w_p^{\eta,h,N}(\cdot,\omega)\right) \cdot \nabla \varphi_h \underset{\eta \rightarrow 0}{\longrightarrow} \int_{Q_N}  A_{\mathrm{per}}
\left( p +  \nabla w_p^{0,h,N}(\cdot,\omega)\right) \cdot \nabla \varphi_h \ \mbox{almost surely}.
\end{eqnarray}
The same argument gives
\begin{eqnarray}
\label{eq:conv-E1}
E^\eta_1(\omega) &\underset{\eta \rightarrow 0}{\longrightarrow}& \int_{Q_N} \mathrm{div} \Psi (\cdot,\omega) A_{\mathrm{per}}
\left( p +  \nabla w_p^{0,h,N}(\cdot,\omega)\right) \cdot \nabla \varphi_h 
-  \int_{Q_N}  \nabla \Psi(\cdot,\omega)^T   A_{\mathrm{per}}
\left( p +  \nabla w_p^{0,h,N}(\cdot,\omega)\right) \cdot \nabla \varphi_h \nonumber \\
&& \quad  - \int_{Q_N}  A_{\mathrm{per}}
\nabla \Psi(\cdot,\omega) \nabla w_p^{0,h,N}(\cdot,\omega) \cdot \nabla \varphi_h,  \ \ \mbox{almost surely}\\
\label{eq:conv-E2}
E^\eta_2(\omega) &\underset{\eta \rightarrow 0}{\longrightarrow}& \int_{Q_N} \nabla \Psi(\cdot,\omega)^T A_{\mathrm{per}} \nabla \Psi(\cdot,\omega) \nabla w_p^{0,h,N}(\cdot,\omega) \cdot \nabla \varphi_h 
- \int_{Q_N} \mathrm{div} \Psi (\cdot,\omega) \nabla \Psi(\cdot,\omega)^T A_{\mathrm{per}}  \left( \nabla w_p^{0,h,N}(\cdot,\omega)  + p \right) \cdot \nabla \varphi_h \nonumber \\
&& \quad - \int_{Q_N} \mathrm{div} \Psi (\cdot,\omega)  A_{\mathrm{per}} \nabla \Psi(\cdot,\omega) \nabla w_p^{0,h,N}(\cdot,\omega) \cdot \nabla \varphi_h , \ \ \mbox{almost surely}.
\end{eqnarray}

\paragraph{Bounds on $E^\eta_1$, $E^\eta_2$, $E^\eta_3$ and $F_\eta$ for sufficiently small values of $\eta$.}

First we remark that
\begin{eqnarray}
\label{eq:bound-E1}
\left|E^\eta_1(\omega) \right| &\leq& 
\left\|A_\mathrm{per} \right\|_{\left(L^\infty(Q_N) \right)^{d \times d}} \left\|\nabla \varphi_h\right\|_{\left(L^2(Q_N) \right)^{d}} \nonumber  \\
&& \quad \quad \times 
 \left(|Q_N|^{\frac{1}{2}}  +  \left\|\nabla w^{\eta,h,N}_p(\cdot,\omega)  \right\|_{\left(L^2(Q_N) \right)^{d}}\right) 
\left(\left\| \mathrm{div} \Psi \right\|_{C^0(\RR^d,L^\infty(\Omega))} + \left\| \nabla \Psi \right\|_{\left(C^0(\RR^d,L^\infty(\Omega))\right)^{d \times d}} \right) \nonumber \\
&& \quad + \left\|A_\mathrm{per} \right\|_{\left(L^\infty(Q_N) \right)^{d \times d}} \left\| \nabla \Psi \right\|_{\left(C^0(\RR^d,L^\infty(\Omega))\right)^{d \times d}} \left\|\nabla \varphi_h\right\|_{\left(L^2(Q_N) \right)^{d}} \left\|\nabla w^{\eta,h,N}_p(\cdot,\omega)  \right\|_{\left(L^2(Q_N) \right)^{d}}, \nonumber \\
\label{eq:bound-E2}
\left|E^\eta_2(\omega) \right| &\leq& \left\|A_\mathrm{per} \right\|_{\left(L^\infty(Q_N) \right)^{d \times d}} \left\|\nabla \varphi_h\right\|_{\left(L^2(Q_N) \right)^{d}} \nonumber  \\
&& \quad \quad \times   \left(|Q_N|^{\frac{1}{2}}  +  \left\|\nabla w^{\eta,h,N}_p(\cdot,\omega)  \right\|_{\left(L^2(Q_N) \right)^{d}}\right) \left\| \mathrm{div} \Psi \right\|_{C^0(\RR^d,L^\infty(\Omega))} \left\| \nabla \Psi \right\|_{C^0(\RR^d,L^\infty(\Omega))} \nonumber \\
&& \quad + \left\|A_\mathrm{per} \right\|_{\left(L^\infty(Q_N) \right)^{d \times d}} \left\|\nabla \varphi_h\right\|_{\left(L^2(Q_N) \right)^{d}} \nonumber \\
&& \quad \quad \times \left\|\nabla w^{\eta,h,N}_p(\cdot,\omega)  \right\|_{\left(L^2(Q_N) \right)^{d}} \left( 
\left\| \nabla \Psi \right\|_{C^0(\RR^d,L^\infty(\Omega))}^2 + \left\| \nabla \Psi \right\|_{C^0(\RR^d,L^\infty(\Omega))} \left\| \mathrm{div} \Psi \right\|_{C^0(\RR^d,L^\infty(\Omega))} \right), \nonumber 
\end{eqnarray}
almost surely. Using~\eqref{eq:Cw-diffeo}, we deduce that there exist two constants $C_{E_1}$ and $C_{E_2}$ independent of $\eta$, $h$, $N$ and $\omega$ such that,
\begin{equation}
\label{eq:CE1-CE2}
\left|E^\eta_1(\omega) \right| \leq C_{E_1}\left\|\nabla \varphi_h\right\|_{\left(L^2(Q_N) \right)^{d}} |Q_N|^{\frac{1}{2}} \ \ \mbox{and} \ \  \left|E^\eta_2(\omega) \right| \leq C_{E_2} \left\|\nabla \varphi_h\right\|_{\left(L^2(Q_N) \right)^{d}}|Q_N|^{\frac{1}{2}}, \ \ \mbox{almost surely}.
\end{equation}
Second, using~\eqref{eq:dl-gamma-sigma}, we observe that each term of $E^\eta_3(\omega)$ is bounded. As an example, the first term of $E^\eta_3(\omega)$ satisfies
\begin{eqnarray}
\left| \int_{Q_N} \sigma_\eta (\cdot,\omega) A_{\mathrm{per}}
\left( p +  \nabla w_p^{\eta,h,N}(\cdot,\omega)\right) \cdot \nabla \varphi_h \right| &\leq& C_\sigma \eta^2 \left\| \nabla \varphi_h \right\|_{\left(L^2(Q_N)\right)^d} \left\| A_\mathrm{per} \right\|_{\left(L^\infty(Q_N)\right)^{d \times d}} \nonumber \\
&& \quad \times \left( \left\| \nabla w_p^{\eta,h,N}(\cdot,\omega)  \right\|_{\left(L^2(Q_N)\right)^d} + |Q_N|^{\frac{1}{2}} \right). \nonumber
\end{eqnarray}
The other terms of $E^\eta_3(\omega)$ can be bounded likewise. Using~\eqref{eq:Cw-diffeo}, we obtain that there exists a constant $C_{E_3}$ independent of $\eta$, $h$, $N$ and $\omega$ such that, for sufficiently small values of $\eta$, we have
\begin{equation}
\label{eq:CEprime2}
|E^\eta_3(\omega)| \leq C_{E_3} \left\| \nabla \varphi_h \right\|_{\left(L^2(Q_N)\right)^d} |Q_N|^{\frac{1}{2}}\eta^2 \ \ \mbox{almost surely}.
\end{equation}
We conclude that 
\begin{equation}
\label{eq:conv-Eprime2}
E^\eta_3(\omega) \underset{\eta \rightarrow 0}{\longrightarrow} 0 \ \ \mbox{almost surely}.
\end{equation}
Bounding from above $F_\eta(\omega)$ is technically more tedious, since it contains many terms. Nevertheless, each of them can be written under one of these two forms
$$ \int_{Q_N} T_\eta (\cdot, \omega) \nabla w_p^{\eta,h,N}(\cdot,\omega) \cdot \nabla \varphi_h \ \ \mbox{or} \ \ \int_{Q_N} T_\eta (\cdot, \omega) \left( \nabla w_p^{\eta,h,N}(\cdot,\omega) + p \right) \cdot \nabla \varphi_h, $$
where $T_\eta(\cdot,\omega) \rightarrow 0 $ strongly in $\left(L^\infty(Q_N) \right)^{d \times d}$ almost surely. As an exemple, let us consider a term corresponding to the choice $T_\eta (\cdot,\omega)=  \eta \mathrm{div} \Psi(\cdot,\omega) A_{\mathrm{per}} \Gamma_\eta (\cdot,\omega)$. Using~\eqref{eq:dl-gamma-sigma} and~\eqref{eq:Cw-diffeo},  we have
\begin{eqnarray}
\left|\int_{Q_N} \eta \mathrm{div} \Psi(\cdot,\omega) A_{\mathrm{per}} \Gamma_\eta (\cdot,\omega) \nabla w_p^{\eta,h,N} \cdot \nabla \varphi_h \right| &\leq&  \eta \left\| \mathrm{div} \Psi(\cdot,\omega)  A_{\mathrm{per}} \Gamma_\eta (\cdot,\omega)\right\|_{\left(L^\infty(Q_N) \right)^{d \times d}} C_{w} |Q_N|^{\frac{1}{2}}\left\|\nabla \varphi_h \right\|_{\left(L^2 (Q_N)\right)^{d}} \nonumber \\
&\leq &\left\|A_{\mathrm{per}} \right\|_{\left(L^\infty(Q_N) \right)^{d \times d}} \left\|\mathrm{div} \Psi(\cdot,\omega) \right\|_{L^\infty(Q_N)} C_\Gamma \eta^3 C_{w} |Q_N|^{\frac{1}{2}}\left\|\nabla \varphi_h \right\|_{\left(L^2 (Q_N)\right)^{d}}. \nonumber
\end{eqnarray}
This term tends to $0$ as $\eta \rightarrow 0$. Carrying over the same analysis for each term of $F_\eta(\omega)$, we conclude that 
there exists a \emph{deterministic} constant $C_{F}$ independent of $h$, $N$ and $\eta$ such that, for sufficiently small values of $\eta$, we have
\begin{equation}
\label{eq:CEeta}
|F_\eta(\omega)| \leq C_{F} \eta^3 \left\|\nabla \varphi_h\right\|_{\left(L^2(Q_N) \right)^{d}} |Q_N|^{\frac{1}{2}}  \ \ \mbox{almost surely}. \\ 
\end{equation}
Of course, this implies
\begin{equation}
\label{eq:conv-Eeta}
F_\eta(\omega)  \underset{\eta \rightarrow 0}{\longrightarrow} 0 \ \ \mbox{almost surely}.
\end{equation}

\paragraph{Zero order equation.} We now return to equation~\eqref{eq:PC-dev-diffeo} and pass to the limit $\eta \rightarrow 0$.
Using~\eqref{eq:CE1-CE2}, it is obvious that
\begin{equation}
\label{eq:conv-E1-bis}
\eta E^\eta_1(\omega) \underset{\eta \rightarrow 0}{\longrightarrow} 0, \ \ \ \ \ \eta^2 E^\eta_2(\omega) \underset{\eta \rightarrow 0}{\longrightarrow} 0 \ \ \mbox{almost surely}. 
\end{equation}
Collecting~\eqref{eq:conv-E0},~\eqref{eq:conv-Eprime2},~\eqref{eq:conv-Eeta} and~\eqref{eq:conv-E1-bis},
we deduce that $w^{0,h,N}_p(\cdot,\omega)$ is a solution to~\eqref{eq:correcteur-0-hN-bis}. We have proved (see the first step of the proof of Proposition~\ref{prop:cas1-dl-hN}) that~\eqref{eq:correcteur-0-hN-bis} has a unique solution which is independent of $N$ and is equal to $w^{0,h}_p$ the unique solution $w^{0,h}_p$ to~\eqref{eq:correcteur-0-hN}.\\

\paragraph{\underbar{Step 2:}} We introduce the function
\begin{equation}
v^{\eta,h,N}_p(\cdot,\omega) = \frac{w^{\eta,h,N}_p(\cdot,\omega) - w^{0,h}_p}{\eta}. \nonumber
\end{equation}
Our aim is to prove that $\nabla v^{\eta,h,N}_p$ is bounded in $\left(L^2(Q_N)\right)^d$, almost surely. Then, we show that $w^{1,h,N}_p$, the limit of $ v^{\eta,h,N}_p$, is in fact the unique solution to~\eqref{eq:correcteur-1-hN-cas2}.
Substracting equation~\eqref{eq:correcteur-0-hN-bis} from~(\ref{eq:PC-dev-diffeo}), we obtain that $v^{\eta,h,N}_p \in V^h_\mathrm{per} (Q_N)$ is such that, for any $\varphi_h \in V_h^{\mathrm{per}} (Q_N),$
\begin{eqnarray}
\label{eq:P-P0-diffeo}
 0&= & \int_{Q_N}  A_{\mathrm{per}}
 \nabla v_p^{\eta,h,N}(\cdot,\omega) \cdot \nabla \varphi_h  + E^\eta_1(\omega) + \eta E^\eta_2(\omega) + \frac{E^\eta_3(\omega)}{\eta} + \frac{F_\eta(\omega)}{\eta}.
\end{eqnarray}
We are now going to use  the bounds on $E^\eta_1$, $E^\eta_2$, $E^\eta_3$, $F^\eta$ obtained below.
Indeed, choosing $\varphi_h = v^{\eta,h,N}_p(\cdot,\omega)$ in~\eqref{eq:P-P0-diffeo},~\eqref{eq:CE1-CE2},~\eqref{eq:CEprime2}, and~\eqref{eq:CEeta}, we obtain that, for sufficiently small values of $\eta$,
$$\gamma  \left\|\nabla v^{\eta,h,N}_p(\cdot,\omega) \right\|_{\left(L^2 (Q_N)\right)^{d}} \leq \left( C_{E_1} + \eta C_{E_2}  + \eta C_{E_3}  + \eta^2 C_{F}  \right) |Q_N|^{\frac{1}{2}},$$
where $C_{E_1}$, $C_{E_2}$, $C_{E_3}$ and $C_F$ are all independent of $\eta$, $h$, $N$ and $\omega$. 
As a consequence, $v_p^{\eta,h,N}(\cdot,\omega)$ is bounded in $(L^2(Q_N))^d$ almost surely. Thus, there exists a deterministic constant $C_v$ independent of $h$, $N$ and $\eta$, such that, for sufficiently small values of $\eta$,
\begin{eqnarray}
\label{eq:C-v-tilde}
\left\|\nabla v^{\eta,h,N}_p(\cdot,\omega) \right\|_{\left(L^2 (Q_N)\right)^{d}} &\leq &C_{v}  |Q_N|^{\frac{1}{2}} \ \ \mbox{almost surely}.
\end{eqnarray}
Thus, there exists a function $w^{1,h,N}_p(\cdot,\omega) \in V^h_\mathrm{per} (Q_N)$ such that, up to extracting a subsequence,
\begin{equation}
\label{eq:conv-v}
 \nabla v^{\eta,h,N}_p (\cdot,\omega) \mathop{\longrightarrow}_{\eta \rightarrow 0} \nabla w^{1,h,N}_p (\cdot,\omega) \ \mbox{in} \left(L^2(Q_N)\right)^{d} \ \mbox{almost surely.} 
\end{equation}
Using~\eqref{eq:conv-E1},~\eqref{eq:CE1-CE2},~\eqref{eq:CEprime2},~\eqref{eq:CEeta} and~\eqref{eq:conv-v}, we can pass to the limit $\eta \rightarrow 0$ in~\eqref{eq:P-P0-diffeo}. We conclude that $w^{1,h,N}_p  $ is the unique function such that, for any $\varphi_h \in V^h_\mathrm{per}(Q_N)$
\begin{eqnarray}
& &\int_{Q_N}  A_{\mathrm{per}}
 \nabla w_p^{1,h,N}(\cdot,\omega) \cdot \nabla \varphi_h - \int_{Q_N}  A_{\mathrm{per}}
\nabla \Psi(\cdot,\omega) \nabla w_p^{0,h} \cdot \nabla \varphi_h  \nonumber \\
&+&\int_{Q_N} \left( \mathrm{div} \Psi (\cdot,\omega) \mathbb{I}_d - \nabla \Psi(\cdot,\omega)^T \right) A_{\mathrm{per}}
\left( p +  \nabla w_p^{0,h}\right) \cdot \nabla \varphi_h =0 . \nonumber
\end{eqnarray}
We recover the problem~\eqref{eq:correcteur-1-hN-cas2}, the solution of which is unique.
\paragraph{\underbar{Step 3:}}  We define 
\begin{equation}
z^{\eta,h,N}_p(\cdot,\omega) = \frac{w_p^{\eta,h,N}(\cdot,\omega) - w^{0,h}_p - \eta w^{1,h,N}_p(\cdot,\omega)}{\eta^2}. \nonumber
\end{equation}
Using~\eqref{eq:correcteur-0-hN-bis},~\eqref{eq:correcteur-1-hN-cas2},~(\ref{eq:PC-dev-diffeo}) and~\eqref{eq:conv-v} we obtain that $z^{\eta,h,N}_p$ is such that, for any $\varphi_h \in V_h^{\mathrm{per}} (Q_N),$
\begin{eqnarray}
\label{eq:P-P0-P1-diffeo}
 0&= & \int_{Q_N}  A_{\mathrm{per}}
 \nabla z_p^{\eta,h,N}(\cdot,\omega) \cdot \nabla \varphi_h + E^\eta_2(\omega) + \frac{E^\eta_3(\omega)}{\eta^2} + \frac{F_\eta(\omega)}{\eta^2}.
\end{eqnarray}
Choosing $\varphi_h = z^{\eta,h,N}_p(\cdot,\omega)$ in~\eqref{eq:P-P0-P1-diffeo}, and using~\eqref{eq:CE1-CE2},~\eqref{eq:CEprime2}, and~\eqref{eq:CEeta}, we obtain that, for sufficiently small values of $\eta$,
\begin{equation}
\gamma \left\|\nabla z^{\eta,h,N}_p(\cdot,\omega) \right\|_{\left(L^2 (Q_N)\right)^{d}} \leq \left(C_{E_2}  + C_{E_3} +  \eta {C_{F}}\right)  |Q_N|^{\frac{1}{2}}. \nonumber
\end{equation}
%This shows that $t_p^{\eta,h,N}(\cdot,\omega)$ is bounded in $(L^2(Q_N))^d$ almost surely. More precisely
We thus have
\begin{equation}
\label{eq:C-t-tilde}
\left\|\nabla z^{\eta,h,N}_p(\cdot,\omega) \right\|_{\left(L^2 (Q_N)\right)^{d}} \leq C_{z} |Q_N|^{\frac{1}{2}},
\end{equation}
where $C_z$
is independent of  $h$, $N$, $\eta$, and $\omega$. Thus, inequality~\eqref{eq:dl-w-hN-cas2} is proved.

\paragraph{\underbar{Step 4:}} Our last step aims at bounding the second order error of the expansion of $A^{\star,h,N}_\eta$. We insert the expansions~\eqref{eq:dl-gradphi} and~\eqref{eq:dl-detphi} into~\eqref{eq:a-hN-cas2}. We obtain that
\begin{eqnarray}
\left[A^{\star,h,N}_\eta \right]_{ij}(\omega)&=&
\frac{1}{|Q_N|}\int_{Q_N} \left(e_i + \nabla w^{\eta,h,N}_{e_i}(\cdot,\omega)\right)^T A_{\mathrm{per}} e_j \nonumber \\
&+& \eta \left(
\begin{array}{ccl}
 &-&\displaystyle{\left( \frac{1}{|Q_N|}\int_{Q_N}\left(e_i + \nabla w^{\eta,h,N}_{e_i}(\cdot,\omega) \right)^T A_{\mathrm{per}} e_j\right) \frac{1}{|Q_N|}\int_{Q_N}  \mathrm{div} \Psi(\cdot,\omega)} \\
&+& \displaystyle{\frac{1}{|Q_N|}\int_{Q_N} \mathrm{div} \Psi (\cdot,\omega) \left(e_i + \nabla w^{\eta,h,N}_{e_i}(\cdot,\omega) \right)^T A_{\mathrm{per}} e_j }\\
&- &\displaystyle{\frac{1}{|Q_N|}\int_{Q_N} \left( \nabla \Psi (\cdot,\omega) \nabla w^{\eta,h,N}_{e_i}(\cdot,\omega) \right)^T A_{\mathrm{per}} e_j}
\end{array}
\right) \nonumber \\
&+& \delta^{h,N}_\eta (\omega), \nonumber
\end{eqnarray}
where $\delta^{h,N}_\eta(\omega)$ is a sum of terms of the forms
$$ s_\eta(\omega) \frac{1}{|Q_N|}\int_{Q_N} T_\eta (\cdot, \omega) \left( \nabla w_{e_i}^{\eta,h,N}(\cdot,\omega) \right)^T A_\mathrm{per}  \cdot e_j \ \ \mbox{and} \ \ s_\eta(\omega) \frac{1}{|Q_N|} \int_{Q_N} T_\eta (\cdot, \omega) \left( \nabla w_{e_i}^{\eta,h,N}(\cdot,\omega) + e_i \right)^T A_\mathrm{per}  \cdot e_j,$$
with $s_\eta(\omega) \in \RR $ and $T_\eta \in C^0\left(\RR^d,L^\infty(\Omega) \right)^{d \times d}$. For each of these $\left(s_\eta,T_\eta \right)$, there exists a constant $C_{s,T}$ 
 independent of  $h$, $N$, $\eta$, and $\omega$ such that, for sufficiently small values of $\eta$,
$$ \left\|s_\eta \right\|_{L^\infty(\Omega)}  \left\|T_\eta \right\|_{\left( C^0 \left( Q_N, L^\infty(\Omega) \right) \right)^{d \times d}} \leq C_{s,T} \eta^2.$$
Therefore, for each term of $\delta_\eta^{h,N}(\omega)$ we can write
\begin{eqnarray}
\left| s_\eta(\omega) \frac{1}{|Q_N|} \int_{Q_N} T_\eta (\cdot, \omega) \left( \nabla w_{e_i}^{\eta,h,N}(\cdot,\omega) \right)^T A_\mathrm{per}  \cdot e_j \right| &\leq& \frac{1}{|Q_N|} \left\|s_\eta \right\|_{L^\infty(\Omega)}  \left\|T_\eta \right\|_{\left( C^0 \left( \RR^d, L^\infty(Q_N) \right) \right)^{d \times d}} \nonumber \\
&\times& \left\|\nabla w^{\eta,h,N}_p(\cdot,\omega) \right\|_{\left(L^2 (Q_N)\right)^{d}}\left\|A_\mathrm{per}\right\|_{\left(L^\infty(Q_N)\right)^{d \times d}}  |Q_N|^{\frac{1}{2}} \nonumber \\
& \leq & C_{s,T} C_{w} \left\|A_\mathrm{per}\right\|_{\left(L^\infty(Q_N)\right)^{d \times d}} \eta^2. \nonumber
\end{eqnarray}
As a consequence, there exists a constant $C_\delta$ independent of $h$, $N$, $\eta$ and $\omega$, such that, for sufficiently small values of $\eta$, 
\begin{equation}
\label{eq:C-delta}
|\delta^{h,N}_\eta(\omega)| \leq C_\delta \eta^2 \ \ \mbox{almost surely}.
\end{equation}
Next, we compute
\begin{eqnarray}
 \eta^{-2} \left(
A^{\star,h,N}_{\eta}(\omega) - A_\mathrm{per}^{\star,h} - \eta A^{\star,h,N}_1(\omega) \right)&=& 
\frac{1}{|Q_N|} \int_{Q_N} \left(e_i + \nabla z^{\eta,h,N}_{e_i}(\cdot,\omega)\right)^T A_{\mathrm{per}} e_j \nonumber \\
&+& \left(
\begin{array}{ccl}
 &-&\displaystyle{\left( \frac{1}{|Q_N|}\int_{Q_N}\left(\nabla v^{\eta,h,N}_{e_i}(\cdot,\omega) \right)^T A_{\mathrm{per}} e_j\right) \frac{1}{|Q_N|}\int_{Q_N}  \mathrm{div} \Psi(\cdot,\omega)} \\
&+& \displaystyle{\frac{1}{|Q_N|}\int_{Q_N} \mathrm{div} \Psi (\cdot,\omega) \left(\nabla v^{\eta,h,N}_{e_i}(\cdot,\omega) \right)^T A_{\mathrm{per}} e_j }\\
&- &\displaystyle{\frac{1}{|Q_N|}\int_{Q_N} \left( \nabla \Psi (\cdot,\omega) \nabla v^{\eta,h,N}_{e_i}(\cdot,\omega) \right)^T A_{\mathrm{per}} e_j}
\end{array}
\right) \nonumber \\
&+& \frac{\delta^{h,N}_\eta (\omega)}{\eta^2} \nonumber.
\end{eqnarray}
Using bounds~\eqref{eq:C-v-tilde},~\eqref{eq:C-t-tilde} and~\eqref{eq:C-delta}, we obtain~\eqref{eq:dl-a-hN-cas2}. This concludes the proof.

\end{preuveter}

\begin{remarque}
As in Section 2, we can pass to the limit $h \rightarrow 0$ in the bounds of Proposition 3.1 and prove that~\eqref{eq:dl-a-hN-cas2} extends to the case when only truncation is taken into account.
Letting $N$ go to infinity is more difficult. Indeed, we need to show that $  \underset{N \rightarrow \infty}{\lim}\left( \underset{h \rightarrow 0}{\lim} A^{\star,h,N}_\eta (\omega) \right)= A^{\star}_\eta, \ \mbox{almost surely}$. Adapting the arguments of~\cite{BourgeatPiatnitski} to Model 2 will be the subject of a future publication.
\end{remarque}

\paragraph{Acknowledgments}
The author thanks Fr\'ed\'eric Legoll for helpful discussions. This work is partially suppported by ONR under contract Grant 00014-09-1-0470.

\pagebreak

\end{document}